\input ./style/arxiv-general.cfg
\documentclass[aap,MSNbibl,secthm,seceqn,dvips]{arximspdf}
\makeatletter
   \@ifpackageloaded{graphicx}{}{\usepackage{graphicx}}
\makeatother

\doi{10.1214/15-AAP1096}
\volume{26}
\issue{1}
\pubyear{2016}
\firstpage{425}
\lastpage{461}
\docsubty{FLA}

\makeatletter
\newcommand{\eqref}[1]{(\ref{#1})}
\newcommand{\rrvert}{\vert}
\newcommand{\rrVert}{\Vert}
\newcommand{\llvert}{\vert}
\newcommand{\llVert}{\Vert}
\newtheorem{prop}{Proposition}[section]
\newtheorem{lemma}{Lemma}[section]
\newtheorem{meta}{Meta Theorem}[section]
\newproclaim{defn}{Definition}[section]
\newproclaim{rmk}{Remark}[section]
\newproclaim{example}{Example}[section]
\makeatother

\begin{document}
\begin{frontmatter}

\title{Rough path recursions and diffusion approximations}
\runtitle{Rough path recursions}

\begin{aug}
\author[A]{\fnms{David}~\snm{Kelly}\corref{}\thanksref{T1}\ead[label=e1]{dtbkelly@gmail.com}\ead[label=u1,url]{http://www.dtbkelly.com}}
\runauthor{D. Kelly}
\affiliation{University of North Carolina, Chapel Hill}
\address[A]{Courant Institute of\\
\quad Mathematical Sciences\\
New York University\\
New York, New York 10012\\
\printead{e1}\\
\printead{u1}}
\end{aug}
\thankstext{T1}{Supported by ONR Grant N00014-12-1-0257.}
%
%
\received{\smonth{2} \syear{2014}}
%
\revised{\smonth{12} \syear{2014}}

\begin{abstract}
In this article, we consider diffusion approximations for a general
class of stochastic recursions. Such recursions arise as models for
population growth, genetics, financial securities, multiplicative time
series, numerical schemes and MCMC algorithms. We make no particular
probabilistic assumptions on the type of noise appearing in these
recursions. Thus, our technique is well suited to recursions where the
noise sequence is not a semi-martingale, even though the limiting noise
may be. Our main theorem assumes a weak limit theorem on the noise
process appearing in the random recursions and lifts it to diffusion
approximation for the recursion itself. To achieve this, we approximate
the recursion (pathwise) by the solution to a stochastic equation
driven by piecewise smooth paths; this can be thought of as a pathwise
version of backward error analysis for SDEs. We then identify the limit
of this stochastic equation, and hence the original recursion, using
tools from rough path theory. We provide several examples of diffusion
approximations, both new and old, to illustrate this technique.
\end{abstract}

\end{frontmatter}

\section{Introduction}\label{sec1}

In this article, we consider the limiting behaviour for a class of
stochastic recursions. These recursions are natural approximations of
continuous time stochastic equations. They arise as models for fast,
discretely evolving random phenomena \cite{guess77,watterson64,norman74,vervaat79,gottwald13} and also as
numerical discretizations of continuous stochastic equations \cite{kloden92}. The class is similar to the rough path schemes of
\cite{davie07} (see also \cite{frizhairer13}, Section~8.5), but more general
in the sense that the noise driving the recursion is not required to be
a rough path, but may be an approximation (or discretization) of a
rough path.

Let $V \dvtx  \mathbb{R}^e \to\mathbb{R}^{e\times d}$ and $\mathbb{V}=
(\mathbb{V}_1,\ldots,\mathbb{V}
_e)$ where $\mathbb{V}_\kappa\dvtx  \mathbb{R}^e \to\mathbb{R}^{d\times
d}$ is defined
by $\mathbb{V}_\kappa^{\alpha\beta}(\cdot) = \sum_\gamma\partial
^\gamma V_\kappa
^{\beta}(\cdot) V^\alpha_\gamma(\cdot)$ for $\alpha,\beta= 1,
\ldots, d$
and $V = (V_\kappa^\beta)$ for $\kappa= 1, \ldots, e$, $\beta= 1,
\ldots,
d$. For each $n\geq0$, define $Y^n_k \in\mathbb{R}^e$ by the recursion
%
\begin{equation}
\label{erecu} Y^n_{k+1} = Y^n_k +
V \bigl(Y^n_k \bigr)\xi^n_{k} +
\mathbb{V} \bigl(Y^n_k \bigr) \dvtx \Xi
^n_{k} + \mbox{error},
\end{equation}
where $\xi^n_{k} \in\mathbb{R}^d$, $\Xi^n_{k} \in\mathbb
{R}^{d\times d}$ are
noise sources and we use the notation $A\dvtx  B = \operatorname{trace} (AB^T) =
\sum_{\alpha, \beta} A^{\alpha\beta} B^{\alpha\beta}$ to denote the
matrix inner product.

Let $\mathcal{P}_n = \{\tau^n_k \dvtx  k=0 ,\ldots, N_n\}$ be a partition of
a finite
time interval $[0,T]$, which gets finer as $n$ tends to infinity. The
vector $\xi^n_{k} $ should be thought of as an approximation of a
random increment
%
\begin{equation}
\label{eincy1} \xi^n_k \approx X \bigl(
\tau^n_{k+1} \bigr) - X \bigl(\tau^n_k
\bigr),
\end{equation}
where $X$ is some given stochastic process (a semi-martingale or
fractional Brownian motion, e.g.). Formally, the symbol
$\approx
$ means that the approximation gets better as $n$ tends to infinity.
Likewise, the matrix $\Xi^n_{k}$ should be thought of as some
approximation of an iterated stochastic integral
%
\begin{equation}
\label{eincy2} \Xi^n_k \approx\int_{\tau^n_k}^{\tau^n_{k+1}}
\bigl(X(s)-X \bigl(\tau ^n_k \bigr) \bigr)\otimes dX(s),
\end{equation}
where $\otimes$ denotes the outer product and where the notion of
stochastic integration (It\^o, Stratonovich or otherwise) is given.

Define the path $Y^n \dvtx [0,T] \to\mathbb{R}^e$ by $Y^n(t) = Y^n_j$ where
$\tau^n_j$ is the largest mesh point in $\mathcal{P}_n$ with $\tau
^n_j\leq t$
[note that we could equally define $Y^n(\cdot)$ by linear\vspace*{1pt}
interpolation, without altering the results of the article]. Our
objective is to show that the path $Y^n(\cdot)$ converges to the
solution of a stochastic differential equation (SDE) driven by $X$ as
$n$ tends to infinity.
%

\begin{rmk}
In the case where $(\xi^n_k,\Xi^n_k)$ are the increments of a rough
path, that is, $\xi^n_k = X(\tau^n_{k+1}) - X(\tau^n_k)$ and $\Xi
^n_k =
\int_{\tau^n_k}^{\tau^n_{k+1}} (X(s)-X(\tau^n_k))\otimes dX(s)$, the
recursions we consider are precisely the rough path schemes defined in
\cite{davie07}. However, we only require that $(\xi^n_k,\Xi^n_k)$ be
\textit{approximations} of rough paths. This means the class of
recursions we consider is much more general than the class of rough
path recursions and includes many natural approximations that do not
fall under \cite{davie07}. This fact will be illustrated by the
examples below.
\end{rmk}


To see why such diffusion approximations should be possible, it is best
to look at a few examples. The most common variant of \eqref{erecu} is
the ``first-order'' recursion, where $\Xi^n=0$, so that
%
\begin{equation}
\label{erec} Y^n_{k+1} = Y^n_k +
V \bigl(Y^n_k \bigr)\xi^n_{k} +
\mbox{error}.
\end{equation}
This resembles an Euler scheme with approximated noise $\xi
^n_k \approx X(\tau^n_{k+1})-X(\tau^n_{k})$. Hence, it is reasonable to believe
that there should be a diffusion approximation $Y^n \Rightarrow Y$
(where $\Rightarrow$ denotes weak convergence of random variables),
where $Y$ satisfies the SDE
%
\[
dY = V(Y)\star dX,
\]
and $\star \,dX$ denotes some method of stochastic integration (e.g., It\^{o},\break 
Stratonovich or otherwise). It turns out that the choice of
approximating sequence\vspace*{1pt} $\xi^n_k$ of the increment $X(\tau^n_{k+1}) -
X(\tau^n_k)$ has a huge influence as to what type of stochastic
integral arises in the limit.

We now explore this idea with a few examples. The first four examples
are first-order recursions as in \eqref{erec} and the final two are
higher order recursions, as in~\eqref{erecu}.

\begin{example}[(Euler scheme)]
Suppose that $B$ is a $d$-dimensional Brownian motion, let $\xi^n_k =
B((k+1)/n)-B(k/n)$ and define the partition $\mathcal{P}_n$ with $\tau
^n_k =
k/n$. Then clearly $Y^n$ defines the usual Euler--Maruyama scheme on
the time window $[0,1]$. It is well known that $Y^n\Rightarrow Y$ where
$Y$ satisfies the It\^{o} SDE
%
\[
dY = V(Y)\,dB.
\]
This creates a feeling that any Euler looking scheme, like \eqref{erec}, should produce It\^o integrals. As we shall see in the next
few examples, when some \textit{correlation} is introduced to the random
variables $\xi^n_k$, this is certainly not the case.
\end{example}

Less trivial recursions of the form \eqref{erec} have shown up in the
areas of population genetics \cite{watterson64,guess77}, econometric
models \cite{vervaat79}, psychological learning models \cite{norman74},
nonlinear time series models \cite{fan03} and MCMC algorithms \cite
{roberts97}, to name but a few. Here, we will list the example from
\cite{guess77}; our analysis follows that performed in~\cite{kushner81}.

\begin{example}[(Population and genetics models)]
In \cite{guess77}, the authors consider the stochastic difference equation
%
\begin{equation}
\label{eguess} Y^n_{k+1} = f \bigl(S^n_k
\bigr) + \exp \bigl(g \bigl(S^n_k \bigr)
\bigr)Y^n_k,
\end{equation}
where $f(0)=g(0) = 0$ and $\{S^n_k\}_{k\geq 0}$ is a stationary sequence of
random variables with $\mathbf{E}S^n_k = \mu/n$, $\operatorname
{var}(S^n_k) = \sigma^2 /
n$, $\operatorname{cov}(S^n_k,S^n_0) = \sigma^2 r_k / n$ and with
mixing assumptions
on the centered sequences $(S^n_k - \mathbf{E}S^n_k)$ and $((S^n_k)^2
- \mathbf{E}
(S^n_k)^2)$. This recursion arises naturally in models for population
growth and also gene selection, where the environment is evolving in a
random way.

Since the equation \eqref{eguess} is linear, the solution can be
written down explicitly. As a consequence, it is easy to directly
identify the limiting behaviour of each term appearing in the solution,
for instance with the help of Prokhorov's theorem. Alternatively, we
can incorporate the problem into the scope of this article by making~\eqref{eguess} look more like the recursion \eqref{erec}. We first write
\begin{eqnarray*}
Y^n_{k+1} &=& Y^n_k + \mathbf{E}f
\bigl(S^n_k\bigr) + \mathbf{E}\bigl(\exp \bigl(g
\bigl(S^n_k\bigr)\bigr)-1\bigr) Y^n_k
\\
&&{}+ \bigl({f\bigl(S^n_k\bigr) - \mathbf{E}f
\bigl(S^n_k\bigr)} \bigr) + \bigl({\exp \bigl(g
\bigl(S^n_k\bigr)\bigr) - \mathbf{E} \exp\bigl(g
\bigl(S^n_k\bigr)\bigr)} \bigr)Y^n_k.
\end{eqnarray*}
Now if we replace $g, f$ and $\exp$ by their second-order Taylor
expansion, we obtain
\begin{eqnarray*}
Y^n_{k+1} &=& Y^n_k +
n^{-1} \biggl(f_s(0)\mu+ \frac{1}{2}
f_{ss}(0)\sigma^2 + \mu g_s(0)Y^n_k
+ \frac{1}{2}\bigl(g_{ss}(0)+g^2_s(0)
\bigr)\sigma^2 Y^n_k \biggr)
\\
&&{}+ n^{-1/2} \bigl( f_s(0) n^{1/2}
\bigl(S^n_k - \mathbf{E}S^n_k
\bigr) + g_s(0)n^{1/2}\bigl(S^n_k -
\mathbf{E}S^n_k\bigr)Y^n_k
\bigr)
\\
&&{}+\frac{n^{-1}}{2} \bigl( f_{ss}(0)n\bigl(\bigl(S^n_k
\bigr)^2 - \mathbf {E}\bigl(S^n_k
\bigr)^2\bigr) \\
&&\hspace*{36pt}{}+ \bigl(g_{ss}(0)+ g_s^2(0)
\bigr)n\bigl(\bigl(S^n_k\bigr)^2 - \mathbf{E}
\bigl(S^n_k\bigr)^2\bigr)
Y^n_k \bigr).
\end{eqnarray*}
Thus, if we set
\begin{eqnarray*}
V^1(y) &=& f_s(0)\mu+ \tfrac{1}{2}
f_{ss}(0)\sigma^2 + \mu g_s(0)Y^n_k
+ \tfrac{1}{2}\bigl(g_{ss}(0)+g^2_s(0)
\bigr)\sigma^2 y,
\\
V^2(y) &=& f_s(0) + g_s(0)y,\qquad
V^3(y) = \tfrac{1}{2} \bigl(f_{ss}(0) +
\bigl(g_{ss}(0)+ g_s^2(0)\bigr)y \bigr)
\end{eqnarray*}
and
\begin{eqnarray*}
\xi^{n,1}_k &=& n^{-1}, \qquad \xi^{n,2}_k
= n^{-1/2} \bigl( n^{1/2} \bigl(S^n_k -
\mathbf{E}S^n_k\bigr) \bigr),
\\
\xi^{n,3}_k &=& n^{-1} \bigl( n \bigl(
\bigl(S^n_k\bigr)^2 - \mathbf{E}
\bigl(S^n_k\bigr)^2\bigr) \bigr),
\end{eqnarray*}
then $Y^n$ satisfies the recursion
%
\[
Y^n_{k+1} = Y^n_k +
V^1(Y_k)\xi^{n,1}_k +
V^2(Y_k)\xi^{n,2}_k +
V^3(Y_k)\xi^{n,3}_k + \mbox{error}.
\]
Moreover, due to the assumptions on $S^n_k - \mathbf{E}S^n_k$ and $(S^n_k)^2
- \mathbf{E}(S^n_k)^2$, the functional central limit theorem for stationary
mixing sequences implies that
%
\[
\Biggl( \sum_{i=0}^{\lfloor n\cdot \rfloor-1}
\xi^{n,2}_i,\sum_{i=0}^{\lfloor n\cdot \rfloor-1}
\xi^{n,3}_i \Biggr) \Rightarrow (W_2,W_3),
\]
where $W_2,W_3$ are Brownian motions with a computable covariance
structure. Thus, we should expect a diffusion limit of the form
%
\[
dY = V^1(Y)\,dt + \bigl(V^2,V^3 \bigr) (Y)
\star(dW_2,dW_3).
\]
By writing down the solution explicitly, it is shown in \cite{guess77}
that this is indeed the case and the method of integration involves a
correction term that is neither It\^o nor Stratonovich.
\end{example}

In nonlinear scenarios, more sophisticated machinery is required \cite{kushner81}, but this still entails quite heavy and
often nonrealistic
mixing assumptions on the stationary sequence. The framework of
martingale problems \cite{varadhan06} has proved quite suitable for
this analysis \cite{kushner81}. In \cite{kurtz91}, the authors
beautifully address the case where the noise is a semi-martingale
sequence, using the idea of a \textit{good} sequence of semi-martingales.
The following example is taken directly from \cite{protter92,kurtz91}.

\begin{example}[(Discrete time asset pricing)]\label{eggood}
Let $r^n_k$ denote the periodic rate of return for a security with
value $S^n_k$. It follows that
%
\[
S^n_{k+1} = S^n_k +
S^n_k r^n_k.
\]
The authors consider the case where $r^n_k$ is a semi-martingale
difference sequence defined in such a way that, if $M^n$ denotes the
partial sum process
%
\[
M^n(t) = \sum_{i=0}^{\lfloor nt \rfloor-1}
r^n_k
\]
then $M^n \Rightarrow M$ where $M$ is a semi-martingale. It is natural
to expect a diffusion approximation
%
\[
dS = S\star dM
\]
with some undetermined method of integration $\star\, dM$. In \cite
{kurtz91}, the authors provide a natural condition on $M^n$ that
ensures this method of integration is It\^o type, which is clearly the
most natural. Sequences $M^n$ that satisfy this condition are called
\textit{good} semi-martingales. Thus, if $M^n$ is good, then
%
\[
dS = S\,dM,
\]
where the integral is of It\^o type. The authors also permit for a
class of semi-martingales which are a reasonable perturbation of a good
semi-martingale. For instance, suppose that $M^n = {\tilde{M}}^n + A^n$,
where ${\tilde{M}}^n$ is a good sequence of semi-martingales and $A^n$
is a
sequence of semi-martingales with $A^n \Rightarrow0$ as $n\to\infty$
(hence ${\tilde{M}}^n \Rightarrow M$). Now define $H^n(t) = \int_0^t A^n(s)
\,dA^n(s)$ where the integral is of It\^o type and $K^n(t) = [{\tilde{M}}
^n,A^n]_t$ where $[\cdot,\cdot]$ denotes quadratic covariation and
suppose that $(M^n,A^n,H^n,K^n) \Rightarrow(M,0,H,K)$ as $n\to\infty
$. Then $S^n \Rightarrow S$ where $S$ satisfies the It\^o SDE
%
\[
dS = S \,dM + S\,d(H-K).
\]
So formally speaking, we have $\star \,dM = dM + d(H-K)$. Thus, two
equally reasonable approximations of $M$ can yield two vastly different
limiting diffusions. This class of perturbed semi-martingales is
comprehensive enough to cover virtually every diffusion approximation
where the recursion is driven by a semi-martingale sequence.
\end{example}

The next example is a rather important one, which unfortunately does
not fit into the classes of diffusion approximations already studied in
the literature. Understanding the diffusion approximation for this
example is one of the main motivations of this paper.

\begin{example}[(Fast--slow systems)]\label{egfastslow}
Let $T \dvtx  \Lambda\to\Lambda$ describe a chaotic dynamical system with
invariant ergodic measure $\mu$. Define the fast--slow system
%
\[
Y^n_{k+1} = Y^n_k +
n^{-1/2} h \bigl(Y^n_k,T^k \omega
\bigr) + n^{-1}f \bigl(Y^n_k,T^k
\omega \bigr),
\]
where $\omega\in\Lambda$ and $h,f \dvtx  \mathbb{T}^e \times\Lambda\to
\mathbb{T}
^e$ where $\mathbb{T}$ denotes the torus $[0,2\pi)$ and $h$ satisfies the
centering condition $\int h(x,y) \mu(dy) = 0$. If we assume that
$\omega
$ is a random variable with law $\mu$, then the path $Y^n(\cdot) =
Y^n_{\lfloor n \cdot \rfloor}$ becomes a random variable on c\`adl\`ag space.
Note that the assumption $\omega\sim\mu$ simply means that the
chaotic dynamical system is started in stationarity.

Fast--slow systems of this type have been considered in \cite{dolgopyat04,gottwald13,liverani14} and are fundamental to the
understanding of naturally occurring physical systems with separated
time scales \cite{majda99}. Previous attempts at diffusion
approximations typically involve heavy mixing assumption on the
dynamical system $T$ which are difficult to prove for most reasonable
systems \cite{kushner81}. In \cite{dolgopyat04}, the author develops an
alternative It\^o calculus, but only in the case where $T$ defines a
partially hyperbolic dynamical system. In \cite{gottwald13}, the
authors study the special case where the noise is additive. This allows
them to use path-space continuity properties to lift convergence of the
partial sum process to convergence of $Y^n$.

Let us see how fast--slow systems fit into the recursion framework
\eqref{erec}. Using a Fourier expansion truncated at level $d$, we can
replace $h(x,y)$ with the product $h(x)v(y)$ where $h \dvtx  \mathbb{R}^e
\to
\mathbb{R}^{e\times d}$, $v \dvtx  \Lambda\to\mathbb{R}^d$, $\int
v(y)\mu(dy) = 0$
and similarly replace $f(x,y)$ with $f(x)g(y)$. Hence, we obtain
%
\[
Y^n_{k+1} = Y^n_k +
n^{-1/2} h \bigl(Y^n_k \bigr)v
\bigl(T^k \omega \bigr) + n^{-1}f \bigl(Y^n_k
\bigr)g \bigl(T^k\omega \bigr).
\]
This clearly satisfies the recursion \eqref{erec} with $V = (h,f)$ and
$\xi^n_k = (n^{-1/2}\times\break  v(T^k\omega), n^{-1}g(T^k\omega))$. The limiting
behaviour of the partial sums
%
\[
W^n(t) = n^{-1/2}\sum_{i=0}^{\lfloor nt \rfloor-1}
v \bigl(T^i \bigr) \quad\mbox{and} \quad S^n(t) = n^{-1}
\sum_{i=0}^{\lfloor nt \rfloor-1} g \bigl(T^i
\bigr)
\]
is well understood under extremely weak conditions on the dynamical
system \cite{young98,melbourne07,melbourne08,melbourne09}. In particular,
%
\[
W^n \Rightarrow W \quad \mbox{and}\quad S^n(t) \to t \bar{g} \qquad\mbox{($\mu$-a.s.)},
\]
where $W$ is a multiple of Brownian motion and $\bar{g}= \int g\, d\mu$.
Thus, we would expect a diffusion approximation of the form
%
\[
dY = h(Y)\star dW + \bar{g}f(Y)\,dt.
\]
In the situations that are already understood, namely partially
hyperbolic dynamical systems \cite{dolgopyat04} or additive noise
\cite{gottwald13}, the limiting stochastic integral shown to be neither It\^
o nor Stratonovich type. Thus, the interpretation of the integral in a
more general setting is an important problem.
\end{example}

The more general family of recursions defined in \eqref{erecu} (with
$\Xi^n_k \neq0$) arise when using a second-order approximation.
Naturally, it is easy to find examples from numerical analysis.

\begin{example}[(Semi-implicit numerical schemes)]\label{egtheta}
Let $X$ be some stochastic process (e.g., fractional Brownian
motion) and introduce the shorthand $X(s,t) = X(t) - X(s)$. Suppose we
approximate a stochastic equation using a semi-implicit method of
integration, for instance,
%
\[
Y^n_{k+1} = Y^n_k +
\tfrac{1}{2} \bigl( V \bigl(Y^n_k \bigr) + V
\bigl(Y^n_{k+1} \bigr) \bigr) X \bigl(\tau^n_k,
\tau^n_{k+1} \bigr).
\]
It is easy to show that $Y^n$ satisfies \eqref{erecu} with $\xi^n_k =
X(\tau^n_k,\tau^n_{k+1})$ and $\Xi^n_k =\frac{1}{2} X(\tau
^n_k,\tau
^n_{k+1})\otimes X(\tau^n_k,\tau^n_{k+1})$. For a simple stochastic
process $X$, like Brownian motion, it is well known that the limit of
this numerical scheme is
%
\[
dY = V(Y) \circ dX,
\]
where the integral is of Stratonovich type. But for more complicated
objects like fractional Brownian motion, it is not so simple \cite{russo93}. Thus, studying recursions of the type \eqref{erecu} can
lead to a better understanding well-posedness for numerical schemes,
that is, whether they are approximating the correct continuous time limit.
\end{example}

\begin{example}[(Sub-diffusion approximations)] \label{egsubdiff}
Despite the article's title, its scope is not restricted to diffusions,
in particular the results also concern sub-diffusions. In \cite{davydov70,taqqu74}, the authors consider partial sum processes of the form
%
\[
X^n(t) = d_n^{-1}\sum
_{k=0}^{\lfloor nt \rfloor-1} \xi_k,
\]
where $\{\xi_k\}_{k\geq0} $ is a stationary dependent sequence of
random variables and $d_n$ is some normalizing constant, such that $X^n
\Rightarrow X$ in the Skorokhod topology, where $X$ is fractional
Brownian motion with some Hurst parameter $H \in(0,1)$ depending on
the correlation structure. With this in mind, it is natural to consider
a recursion
%
\[
Y^n_{k+1} = Y^n_k + V
\bigl(Y^n_k \bigr)\,d_n^{-1}
\xi_k + \mathbb{V} \bigl(Y^n_k \bigr) \dvtx
\Xi^n_k,
\]
where
$\Xi^n_{k}$ is some\vspace*{-1.5pt} approximation of an iterated integral defined using
the sequence $\{ \xi_k\}_{k\geq0}$. For instance, as indicated by
Example~\ref{egtheta}, if $\Xi^n_{k}=\frac{d^{-2}_n}{2}\xi
_k\otimes
\xi_k$ then the above recursion corresponds to a mid-point rule
approximation of a stochastic integral. In particular one would expect
$Y^n \Rightarrow Y$ where
%
\[
Y(t) = Y(0) + \int_0^t V \bigl(Y(s) \bigr)
\circ dX(s)
\]
and where the integral is of \textit{symmetric} type \cite{russo93},
which is the natural limit of the mid-point scheme. Of course, this is
only a guess and it is quite possibly wrong. As we will see, the tools
introduced in this article provide a natural basis for the
investigation of such sub-diffusion approximations. Understanding such
recursions is vastly important and could facilitate for the design of
new methods for simulating stochastic differential equations driven by
fractional Brownian motion.
\end{example}

The technique employed in this article is similar in spirit to that
found in \cite{kurtz91,gottwald13}, in that we will lift an
approximation result for the noise signal into diffusion approximation
for the recursion. However, in our more general scenario, where we do
not assume any particular probabilistic structure on the noise, we
require not just an invariance principle for the noise but also for its
\textit{iterated integral}. More precisely, define the noise signal
%
\[
X^n(t) \stackrel{\mathrm{def}} {=}\sum
_{i=0}^{\lfloor nt \rfloor
-1} \xi^n_i,
\]
which is the natural approximation of the limiting noise signal $X$.
Moreover, define the discrete iterated integral
%
\[
\mathbb{X}^n(t) \stackrel{\mathrm{def}}{=}\sum
_{i=0}^{\lfloor nt
\rfloor-1} \sum_{j=0}^{i-1}
\xi^n_i \otimes\xi^n_j + \sum
_{i=0}^{\lfloor nt \rfloor-1} \Xi^n_i,
\]
which\vspace*{1pt} is the natural approximation of the limiting iterated integral
$\int_0^t X\otimes dX$. In this paper, we shall lift a limit theorem
for the discrete pair $(X^n,\mathbb{X}^n)$ into a diffusion
approximation for
the recursion $Y^n$. In essence, the limiting behaviour of $X^n$ tells
us what type of noise appears in the limiting stochastic integral and
the limiting behaviour of $\mathbb{X}^n$ tells us what type of stochastic
integral we are talking about. Looking back at Example~\ref{egfastslow}, for instance, this suggest that we can interpret the
integral $\star \,dW$, provided we can identify the limit of the discrete
iterated integral
%
\[
\sum_{i=0}^{\lfloor nt \rfloor-1} \sum
_{j=0}^{i-1} v \bigl(T^j \bigr) \otimes v
\bigl(T^i \bigr).
\]
To derive this diffusion approximation technique, we use tools from
\textit{rough path theory} \cite{lyons98}.

\subsection{Diffusion approximations using rough path theory}
For stochastic differential equations driven by piecewise smooth
signals, the relationship between the noise and the solution is
extremely well understood---mostly thanks to {rough path theory}. For
the purpose of exposition, suppose that $X$ is some piecewise smooth
stochastic process and that $Y$ solves the equation
%
\begin{equation}
\label{eintrorde}
Y(t) = Y(0) + \int_0^t V
\bigl(Y(s) \bigr)\,dX(s),
\end{equation}
where the integral is defined in the Riemann--Stieltjes sense. It is
well known that the map $X \mapsto Y$ is not continuous in the sup-norm
topology. The theory of rough path proposes that we can build a
continuous map from the noise to the solution, provided we know a bit
more information about $X$. In particular, suppose that we can define
$\mathbb{X}(t) = \int_0^t X(s)\otimes dX(s)$ where the integral is
again of
Riemann--Stieltjes type. Then one can show that the map $(X,\mathbb{X})
\mapsto Y$ is continuous in a topology called the $\rho_\gamma$
topology (known colloquially as the \textit{rough path topology}). This
topology can be thought of as an extension of the $\gamma$-H\"older
topology, defined on the space of objects similar to the pair
$(X,\mathbb{X}
)$. The objects $(X,\mathbb{X})$ are called rough paths and the metric space
of such objects is called the space of $\gamma$-H\"older rough paths.

This idea clearly has ramifications to the diffusion approximations.
Indeed, suppose that $Y^n$ solves the stochastic equation
%
\[
Y^n(t) = Y^n(0) + \int_0^t
V \bigl(Y^n(s) \bigr)\,dX^n(s)
\]
for some smooth stochastic process $X^n$ and also define the iterated
integral $\mathbb{X}^n(t) = \int_0^t X^n(s) \otimes dX^n(s)$. Since
continuous maps preserve weak convergence, this suggests that a weak
limit theorem for the pair $(X^n,\mathbb{X}^n)$ in the $\rho_\gamma$ topology
can be lifted to a weak limit theorem for $Y^n$. The general procedure
can be summarized by two steps.
\begin{longlist}[2.]
\item[1.] Show that $(X^n,\mathbb{X}^n) \stackrel{\mathrm{f.d.d.}}{\to
}(X,\mathbb{X})$, where $\stackrel{\mathrm{f.d.d.}}{\to}$ denotes
convergence of finite-dimensional distributions.
\item[2.] Show that the sequence is tight in the $\rho_\gamma$ topology.
For instance, one could use a Kolmogorov type argument, by checking
estimates of the form
%
\begin{equation}
\label{eestimates0}
\bigl( \mathbf{E}\bigl|X^n(s,t)\bigr|^q
\bigr)^{1/q} \lesssim|t-s|^\gamma\quad \mbox{and}\quad \bigl( \mathbf{E}\bigl|
\mathbb{X}^n(s,t)\bigr|^{q/2} \bigr)^{2/q} \lesssim
|t-s|^{2\gamma}
\end{equation}
for all $s,t \in[0,T]$, with some suitable $\gamma$ and with $q$
large enough.
\end{longlist}
Since the map: \textit{rough path} $\mapsto$ \textit{solution} is
continuous in the rough path topology, the conclusion from these two
steps is that $Y^n \Rightarrow Y$ where $Y$ is the solution to an SDE
whose form can be determined by the limit $\mathbb{X}$. For instance, suppose
that $X$ were a continuous semi-martingale and that
%
\[
\mathbb{X}(t) = \int_0^t X(r) \circ dX(r) +
\lambda t,
\]
where the above integral is Stratonovich type. Then the limiting
equation can be written
%
\[
dY = V(Y) \circ dX + \lambda\dvtx  \mathbb{V}(Y)\, dt.
\]
This precise idea has proved useful in the areas of stochastic
homogenization \cite{lejay05} and equations driven by random walks
\cite{friz09walk}.

Unfortunately, for the recursion \eqref{erecu} the path $Y^n$ does not
satisfy a stochastic equation in the sense of rough path theory, so we
cannot simply apply the above procedure.

The objective of this article is to overcome this obstacle. It turns
out that the same two step procedure defined above, more or less still
works. All we have to do is replace iterated integrals with their
discrete counterparts and replace step $2$ with the same statement
\textit{up to some resolution}. That is, we need only check the estimates
\eqref{eestimates0} for all $s,t \in\mathcal{P}_n$, which requires no
continuity at all. In checking these discrete estimates, we obtain a
tightness-like result for a discrete version of the H\"older metric,
defined (on c\`adl\`ag paths) by
%
\[
\max_{s\neq t \in\mathcal{P}_n}\frac{|A(t)-A(s)|}{|t-s|^\gamma}.
\]
This is of course \textit{always} finite, since it is a maximum over a
finite set, but the tightness result will tell us something about the
asymptotics.

At the heart of the proof is an approximation theorem (Theorem~\ref{thmmod}), which we believe to be useful in its own right. The theorem
allows us to approximate the recursion \eqref{erecu} with the solution
to a stochastic differential equation driven by piecewise smooth paths.
This approximation theorem can be thought of as a generalization of
\textit{the method of modified equations} for SDEs \cite{zygalakis11}
(otherwise known as \textit{backward error analysis} \cite{debussche12}).
In particular, our approximation theorem has the advantage of being
completely pathwise, without depending on the probabilistic properties
of the stochastic process $X^n$ whatsoever. By approximating $Y^n$ by
the solution to a genuine stochastic equation, we unlock the tools of
rough path theory introduced above.

The outline of the paper is as follows. In Section~\ref{sresults}, we
sketch the main theorem of the paper.
In Section~\ref{sapps}, we list a few applications.
In Section~\ref{srough}, we give a brief introduction to rough path
theory and mention some results that are important to the present
article. In Section~\ref{sdrp}, we rigourously define rough paths
recursions, these are the central objects to the article. In Section~\ref{sproperties}, we derive the properties of rough path recursions
that will be needed for the main theorem. In Section~\ref{sconvergence}, we prove the main theorem of the article, concerning
weak convergence of rough path recursions.

\section{The main results and some applications}\label{sresults}
In this section, we state the main theorem, avoiding the technical
definitions that will be introduced in subsequent sections. In
particular, the main theorem (Theorem~\ref{thmmain}) can be stated and
applied without requiring any knowledge of rough path theory and
similarly for the approximation theorem (Theorem~\ref{thmmod}).

Let $\mathcal{P}_n = \{ \tau^n_j \dvtx  j=0,\ldots, N_n \}$ be a partition
of $[0,T]$
with mesh size $\Delta_n = \max_j |\tau^n_{j+1}-\tau^n_j|$. As stated
above, one should regard $\xi^n_{j} \in\mathbb{R}^d$ as an approximation
of the increment
%
\begin{equation}
\label{einc1} \xi^n_j \approx X \bigl(
\tau^n_{j+1} \bigr) - X \bigl(\tau^n_{j}
\bigr).
\end{equation}
Likewise, one should regard $\Xi^n_{j} \in\mathbb{R}^{d\times d}$ as an
approximation of the iterated integral
%
\begin{equation}
\label{einc2}
\Xi^n_j \approx\int_{\tau^n_j}^{\tau^n_{j+1}}
\bigl(X(s)-X \bigl(\tau ^n_j \bigr) \bigr)\otimes dX(s).
\end{equation}
The only consequence of this analogy is that it influences how we
define the \textit{path} corresponding to the incremental processes.
Indeed, the increments can be anything at all, provided they satisfy
the convergence properties stated in the theorem below. To recap, the
recursions we consider in this article are of the form
%
\begin{equation}
\label{erecursion}
Y^n_{j+1} = Y^n_j
+ V \bigl(Y^n_j \bigr)\xi^n_{j}
+ \mathbb{V} \bigl(Y^n_j \bigr) \dvtx  \Xi
^n_{j} + r^n_{j},
\end{equation}
where $j=0,\ldots,N_n-1$ and $|r^n_j|\lesssim\Delta_n^{\lambda}$ for
some $\lambda>1$ and the implied constant is uniform in $n$.

We now define the \textit{rough step-function} $(X^n,\mathbb{X}^n)$
corresponding to the increments $\xi^n_{j},\Xi^n_{j}$. If $\tau^n_k$ is
the largest grid point in $\mathcal{P}_n$ such that $\tau^n_k \leq t$ then
%
\begin{equation}
\label{epaths1}
X^n(t) = \sum_{j=0}^{k-1}
\xi^n_{j}\quad \mbox{and}\quad \mathbb {X}^n(t) = \sum
_{i=0}^{k-1} \sum
_{j=0}^{i-1} \xi^n_i \otimes
\xi^n_j + \sum_{i=0}^{k-1}
\Xi^n_i.
\end{equation}
We similarly define the incremental paths
%
\begin{equation}
\label{epaths2}
X^n(s,t) = \sum_{j=l}^{k-1}
\xi^n_{j} \quad\mbox{and}\quad \mathbb {X}^n(s,t) =
\sum_{i=l}^{k-1} \sum
_{j=l}^{i-1} \xi^n_i \otimes
\xi^n_j + \sum_{i=l}^{k-1}
\Xi^n_i,
\end{equation}
where $\tau^n_l$ is the largest grid point in $\mathcal{P}_n$ such
that $\tau
^n_l \leq s$. It is easy to check that this is the natural choice,
given the motivation \eqref{einc1} and \eqref{einc2}. The main
theorem is as follows.
%
\begin{thm}\label{thmmain}
Let $Y^n$ satisfy \eqref{erecursion} and let $(X^n,\mathbb{X}^n)$ be
c\`adl\`
ag paths defined by \eqref{epaths1}. Suppose that $(X^n,\mathbb
{X}^n) \stackrel{\mathit{f.d.d.}}{\to}
(X,\mathbb{X})$ where $X$ is a continuous semi-martingale and $\mathbb
{X}$ is of
the form
%
\[
\mathbb{X}(t) = \int_0^t X(r) \otimes\circ\,dX(r) + \nu t,
\]
where the integral is defined in the Stratonovich sense and $\nu\in
\mathbb{R}^{d\times d}$. Suppose that the pair $(X^n,\mathbb{X}^n)$
satisfy the estimates
%
\begin{eqnarray}
\bigl(\mathbf{E}\bigl|X^n \bigl(\tau^n_j,
\tau^n_k \bigr)\bigr|^q \bigr)^{1/q}
&\lesssim &  \bigl|\tau^n_j - \tau^n_k\bigr|^{\gamma}
\quad\mbox{and}
\nonumber
\\[-8pt]
\label{eestimates}
\\[-8pt]
\nonumber
\bigl(\mathbf{E}\bigl|\mathbb {X}^n \bigl(\tau
^n_j,\tau^n_k
\bigr)\bigr|^{q/2} \bigr)^{2/q} & \lesssim & \bigl|\tau^n_j
- \tau ^n_k\bigr|^{2\gamma}
\end{eqnarray}
for all $\tau^n_j,\tau^n_k \in\mathcal{P}_n$ where $q>0$, $\gamma
\in(1/3 +
q^{-1},1/2]$ and the implied constant is uniform in $n$. Then $Y^n
\Rightarrow Y$ in the sup-norm topology, where $Y$ satisfies the SDE
%
\[
dY = V(Y)\circ dX + \nu\dvtx  \mathbb{V}(Y)\, dt.
\]
\end{thm}

\begin{rmk}
Although we only require $q>0$ it is clear from $\gamma\in
(1/3+q^{-1},1/2]$ that we always have $q > 6$.
\end{rmk}

\begin{rmk}\label{rmkqsmall}
If the estimates \eqref{eestimates} hold for $\gamma= 1/2$ and all $q
\geq1$, then the condition $V \in C^{3}$ can be relaxed to $V \in
C^{2+}$. This follows using the standard techniques of $(p,q)$ rough
paths (see \cite{lejay06} and \cite{friz10}, Chapter~12). Additional
details will be given in Remark~\ref{rmkqlarge}.
\end{rmk}

\begin{rmk}\label{rmkdrift1}
The result naturally extends to the case with an additional ``drift''
vector field $W \in C^{1+} (\mathbb{R}^e; \mathbb{R}^e)$
%
\[
Y^n_{j+1} = Y^n_j + V
\bigl(Y^n_j \bigr)\xi^n_{j} +
\mathbb{V} \bigl(Y^n_j \bigr) \dvtx  \Xi ^n_{j}
+ W \bigl(Y^n_j \bigr) \bigl(\tau^n_{j+1}
- \tau^n_j \bigr) + r^n_{j}.
\]
In this setting, the limiting SDE is given by
%
\[
dY = V(Y) \circ dX + \bigl(\nu\dvtx  \mathbb{V}(Y) + W(Y) \bigr)\,dt.
\]
This is a more natural way to treat the problem introduced in Example~\ref{egfastslow}. As with Remark~\ref{rmkqsmall}, this extension is a
standard application of $(p,q)$ rough paths.
\end{rmk}

The next result is not so much a theorem as it is a guide for other
theorems. It applies to situations where the noise driving the limiting
equation is not a semi-martingale, such as the sub-diffusions
encountered in Example~\ref{egsubdiff}.

\begin{meta}\label{thmmeta}
In the same context as above. Suppose that $(X^n,\break \mathbb
{X}^n)\stackrel{\mathit{f.d.d.}}{\to}(X,\mathbb{X}
)$ where $X$ is some continuous stochastic process and
%
\[
\mathbb{X}(t) = \int_0^t X(r) \star dX(r),
\]
where $\star\, dX$ denotes some constructible method of integration.
Suppose moreover that $(X^n,\mathbb{X}^n)$ satisfy the estimates
\eqref{eestimates}. Then $Y^n \Rightarrow Y$ where $Y$ satisfies the
stochastic equation
%
\[
Y(t) = Y(0) + \int_0^t V \bigl(Y(s) \bigr)
\star dX(s).
\]
\end{meta}

Theorem~\ref{thmmain} and Meta Theorem \ref{thmmeta} will be proved
in Section~\ref{sconvergence}. The proof of the meta theorem indicates
what we mean by a ``constructible method of integration''.

Finally, the main tool used to derive the results above is the
approximation theorem, which should be thought of as a pathwise version
of backward error analysis (or the method of modified equations) \cite{debussche12,zygalakis11}. The rate estimate depends on the
\textit{discrete} $\gamma$-H\"older norm $C_{\gamma,n}$ which is the smallest
number such that
%
\begin{equation}
\label{eCn}
\qquad\bigl|X^n \bigl(\tau^n_j,
\tau^n_k \bigr)\bigr| \leq C_{\gamma,n} \bigl|
\tau^n_j - \tau ^n_k\bigr|^{\gamma
}\quad
\mbox{and}\quad \bigl|\mathbb{X}^n \bigl(\tau^n_j,
\tau^n_k \bigr)\bigr|\leq C_{\gamma,n}^2 \bigl|
\tau^n_j - \tau^n_k\bigr|^{2\gamma}
\end{equation}
for all $\tau^n_j,\tau^n_k \in\mathcal{P}_n$. Since this number can be
achieved by taking the maximum over a finite set, it is clear that each
$C_{\gamma,n}$ is finite, regardless of the path $(X^n,\mathbb
{X}^n)$. We
will always need some kind of asymptotic estimate on $C_{\gamma,n}$ to
make use of the approximation theorem.

\begin{thm}\label{thmmod}
Suppose that $Y^n(\cdot)$ is the path defined by the recursion \eqref
{erecursion} and that the pair $(X^n,\mathbb{X}^n)$ defined by \eqref
{epaths1} satisfy the estimates \eqref{eestimates} for some
$q,\gamma
$ as in Theorem~\ref{thmmain}. Then for each $n$ we can find a pair of
piecewise smooth paths $(\tilde{X}^n,\tilde{Z}^n) \dvtx  [0,T] \to\mathbb{R}
^d\times\mathbb{R}^{d\times d}$ such that if $\tilde{Y}^n$ solves
%
\begin{equation}
\label{eXZ}
\tilde{Y}^n(t) = \tilde{Y}^n(0) + \int
_0^t V \bigl(\tilde{Y}^n(s) \bigr)\,d
\tilde {X}^n(s) + \int_0^t
\mathbb{V} \bigl(\tilde{Y}^n(s) \bigr) \dvtx d\tilde{Z}^n(s),
\end{equation}
where the integrals are of Riemann--Stieltjes type, then
%
\begin{equation}
\label{eKn}
\bigl\llVert \tilde{Y}^n - Y^n \bigr
\rrVert _\infty\lesssim K_{\gamma,n} \Delta_n^{3\gamma-1},
\end{equation}
where $\llVert  \cdot \rrVert _\infty$ denotes the sup-norm and
where the constant
$K_{\gamma,n} = 1\wedge C_{\gamma,n}^4$, where $C_{\gamma,n}$ is the
constant defined in \eqref{eCn}.
\end{thm}
The proof of Theorem~\ref{thmmod} is contained in Section~\ref{sproperties}. We will give one simple example to illustrate the idea
behind this approximation theorem.

\begin{example}
Suppose that $B$ is a Brownian motion and that
%
\[
Y^n_{k+1} = Y^n_k + V
\bigl(Y^n_k \bigr) \bigl(B \bigl(\tau^n_{k+1}
\bigr) - B \bigl(\tau ^n_{k} \bigr) \bigr).
\]
It is easy to check that, for almost every Brownian path, the constant
$C_{\gamma,n}$ defined in \eqref{eCn} is bounded uniformly in $n$, for
any $\gamma< 1/2$. It follows that we can find an equation driven by
smooth paths, with solution $\tilde{Y}^n$ such that
%
\[
\bigl\llVert \tilde{Y}^n - Y^n \bigr\rrVert
_\infty\lesssim\Delta _n^{3\gamma-1}
\]
for any $\gamma< 1/2$.
\end{example}

\section{Some applications}\label{sapps}
We will now discuss some potential applications for the tools
introduced above.

\subsection{Random walk recursions}\label{srw}
We start with a quite trivial and well known result, with the sole
intention of illustrating how Theorem~\ref{thmmain} should be used.
A~continuous time version of this example can be found in \cite
{friz09walk}. It should be said that the following can easily be
deduced from either \cite{kurtz91} or \cite{friz09walk}.

Suppose that $\{\xi_i\}_{ i\geq0}$ is an i.i.d. sequence of random
variables taking values in $\mathbb{R}^d$ with $\mathbf{E}\xi_i = 0$
and $\mathbf{E}\xi
_i \otimes\xi_i = D$, with $D \in\mathbb{R}^{d\times d}$. We will
consider the random walk recursion
%
\[
Y^n_{k+1} = Y^n_k +
n^{-1/2}V \bigl(Y^n_k \bigr)\xi_k
\]
with associated partition $\mathcal{P}^n$ with $\tau^n_k = k/n$. If
we define
the path $Y^n(\cdot) = Y^n_{\lfloor n\cdot \rfloor}$ then it is well
known that
$Y^n \Rightarrow Y$ in c\`adl\`ag space (with sup-norm topology), where
$Y$ solves the SDE
%
\[
dY = V(Y)D^{1/2}\,dW,
\]
where $W$ is standard Brownian motion on $\mathbb{R}^d$. The following
lemma illustrates how to prove this using Theorem~\ref{thmmain}.
First, we define the rough step function
%
\[
X^n(t) = n^{-1/2}\sum_{i=0}^{\lfloor nt \rfloor-1}
\xi_i, \qquad \mathbb{X}^n(t) = n^{-1}\sum
_{i=0}^{\lfloor nt \rfloor-1}\sum_{j=0}^{i-1}
\xi _j\otimes\xi_i.
\]

\begin{lemma}\label{lem31}
Suppose that $\mathbf{E}|\xi_0|^q < \infty$ for some $q>6$. Then the pair
$(X^n,\mathbb{X}^n)$ satisfy the assumptions of Theorem~\ref
{thmmain} with
$X = D^{1/2}W$ and $\nu=-\frac{1}{2}D$. In particular $Y^n
\Rightarrow
Y$ where
%
\[
dY = V(Y)D^{1/2}\,dW,
\]
where the integral is of It\^o type.
\end{lemma}

\begin{rmk}
The moment condition on $\mathbf{E}|\xi_0|^q$ is much stronger than required
by more traditional solutions to the problem. This is due to the fact
that the conclusion of the theorem is actually \textit{stronger} than
most traditional versions. In particular, we could actually show that
$Y^n$ converges in (a discrete version) of the rough path topology,
which is much stronger than the sup-norm topology. We will not pursue
such statements in this article.
\end{rmk}

\begin{pf*}{Proof of Lemma~\protect\ref{lem31}}
From Donsker's theorem, we already know that $X^n \Rightarrow X =
D^{1/2}W$. To identify the limit of $\mathbb{X}^n$ we simply write it
as a
stochastic integral. In particular, we see that
%
\[
\mathbb{X}^n(t) = \int_0^t
X^n(s-)\otimes dX^n(s),
\]
where the integral is of left-Riemann type (hence It\^o type). That is,
%
\[
\int_0^t Y(s-)\, dZ(s) = \lim\sum
Y(s_i ) \bigl(Z(s_{i+1}) - Z(s_i) \bigr),
\]
where $\{s_i\}$ is a partition of $[0,t]$ and the limit is taken as the
maximum of $s_{i+1} - s_i$ tends to zero. The theory of good
semi-martingales \cite{kurtz91} provides a class of semi-martingale
sequences for which the limit of a sequence of It\^o integrals is an
It\^o integral.

Since the partial sum process $X^n$ is clearly a martingale with
respect to the filtration generated by the sequence $\{\xi_i\}
_{i=0}^{\lfloor nt \rfloor-1}$, we can appeal to \cite{kurtz91}, Theorem~2.2.
In particular, since the quadratic variation
\[
\bigl[X^n,X^n \bigr]_t = n^{-1}
\sum_{i=0}^{\lfloor nt \rfloor-1} \xi_i \otimes
\xi_i,
\]
we have that $\mathbf{E}[X^n,X^n]_t = D\lfloor nt \rfloor/n$ which is
of course
bounded uniformly in $n$. Thus, $X^n$ is {good} and \cite{kurtz91}, Theorem~2.2,  immediately tells us that
%
\[
\bigl(X^n,\mathbb{X}^n \bigr)  \Rightarrow  (X,\mathbb{X})
\]
in the sup-norm topology, where
%
\[
\mathbb{X}(t) = \int_0^t X(s) \otimes dX(s) =
\int_0^t X(s) \otimes \circ\, dX(s) -
\frac{1}{2}Dt,
\]
where the integrals are of It\^o and Stratonovich type, respectively, and
we have converted between them in the usual way. This is of course
stronger than the finite-dimensional distribution result which we
required, but the tools from \cite{kurtz91} make it quite easy to prove.

All that remains is to obtain the discrete tightness estimates \eqref
{eestimates} with $q>6$ and $\gamma= 1/2$. Since $X^n$ is a
martingale, we can apply the Burkholder--Davis--Gundy (BDG) inequality
\begin{eqnarray*}
\mathbf{E}\bigl|X^n(j/n,k/n)\bigr|^{q} &\lesssim & \mathbf {E}\bigl|
\bigl[X^n,X^n\bigr]_{j/n,k/n}\bigr|^{q/2}
\\
&= & n^{-q/2}\mathbf{E}\Biggl\llvert \sum_{i=j}^{k-1}
\xi_i\otimes\xi_i \Biggr\rrvert ^{q/2} \leq
n^{-q/2} \mathbf{E} \Biggl( \sum_{i=j}^{k-1}
|\xi_i|^2 \Biggr)^{q/2}.
\end{eqnarray*}
By the H\"older inequality, we see that
%
\[
\sum_{i=j}^{k-1} |\xi_i|^2
\leq(k-j)^{1-q/2} \Biggl(\sum_{i=j}^{k-1}
|\xi|^q \Biggr)^{2/q}.
\]
It follows that
\begin{eqnarray*}
\mathbf{E}\bigl|X^n(j/n,k/n)\bigr|^{q} &\lesssim &  n^{-q/2}
(k-j)^{q/2-1} \sum_{i=j}^{k-1}
\mathbf{E}|\xi_i|^q
\\
&= & n^{-q/2} (k-j)^{q/2-1} \sum_{i=j}^{k-1}
\mathbf{E}|\xi_0|^q = \mathbf{E}|\xi_0|^q
(k/n-j/n)^{q/2}.
\end{eqnarray*}
Since $\mathbb{X}^n$ is a stochastic integral (or martingale
transform) it
too is a martingale and hence we can again apply the BDG inequality. A
similar argument yields
%
\[
\mathbf{E}\bigl|\mathbb{X}^n(j/n,k/n)\bigr|^{q/2} \lesssim\mathbf{E}|
\xi _0|^q (k/n-j/n)^{q}.
\]
And since $q>6$, the interval $(1/3+q^{-1},1/2]$ is nonempty, so we do
indeed satisfy the requirements of Theorem~\ref{thmmain}. It follows
that $Y^n \Rightarrow Y$ where
%
\[
dY = V(Y)\circ D^{1/2} \,dW - \tfrac{1}{2} D \dvtx  \mathbb{V}(Y)\, dt
\]
and we obtain the required expression by converting Stratonovich to It\^o.
\end{pf*}

\subsection{Fast--slow systems}
Instead of showing how the tools can be used on each of the examples
given in the \hyperref[sec1]{Introduction}, we concentrate on the fast--slow systems,
since it is the least understood. The tools of this article are applied
to fast--slow systems in a companion paper \cite{kelly14} (see also
\cite
{kelly14b}), to yield new results for fast--slow systems. The dynamical
system theory required is slightly too involved to be included in this
paper, thus we will only sketch the ideas behind the result.

We will restrict our attention to the fast--slow system
%
\[
Y^n_{k+1} = Y^n_k +
n^{-1/2} h \bigl(Y^n_k \bigr)v
\bigl(T^k \omega \bigr),
\]
the general case is treated in \cite{kelly14}. Setting $V = h$, $\xi
^n_k = n^{-1/2} v(T^k)$ and $\Xi^n_k=0$ we see that the rough step
function is defined by
%
\[
X^n(t) = n^{-1/2} \sum_{i=0}^{\lfloor nt \rfloor-1}
v \bigl(T^i \bigr)\quad \mbox{and}\quad \mathbb{X}^n(t) =
n^{-1}\sum_{i=0}^{\lfloor nt \rfloor-1}\sum
_{j=0}^{i-1} v \bigl(T^i \bigr)
\otimes v \bigl(T^j \bigr).
\]
We will also introduce the sigma algebra $\mathcal{M}$ which is
whatever sigma
algebra we chose to go with the measure space $(\Lambda,\mu)$.

\begin{prop}
Under ``sufficient'' mixing conditions on $T$, the pair $(X^n,\mathbb{X}^n)$
satisfy the assumptions of Theorem~\ref{thmmain} with $X = D^{1/2}W$ where
%
\[
D^{\alpha\beta} = \int v^\alpha v^\beta \,d\mu+ \sum
_{j=1}^\infty \biggl(\int v^\alpha
v^\beta \bigl(T^j \bigr)\, d\mu+ \int v^\alpha
\bigl(T^j \bigr) v^\beta \,d\mu \biggr)
\]
and
%
\[
\nu^{\alpha\beta} = -\frac{1}{2} \int v^\alpha v^\beta d
\mu+ \frac
{1}{2} \sum_{j=1}^\infty
\biggl(\int v^\alpha v^\beta \bigl(T^j \bigr) \,d\mu-
\int v^\alpha \bigl(T^j \bigr) v^\beta \,d\mu \biggr).
\]
In particular, $Y^n \Rightarrow Y$ where
%
\[
dY = h(Y)D^{1/2}\circ dW + \nu\dvtx  \mathbb{H}(Y)\, dt,
\]
where $\mathbb{H}$ is defined precisely as $\mathbb{V}$, but in terms
of $h$.
\end{prop}

\begin{pf*}{Sketch of proof}
To identify the limit of the pair $(X^n,\mathbb{X}^n)$, we proceed similarly
to the random walk recursion case, namely identify the limit of $X^n$
and then lift it to $\mathbb{X}^n$. To identify the limit of $X^n$, we will
use a martingale central limit theorem on the \textit{time reversal} of
the partial sum process $X^n$.

By applying the \textit{natural extension} of a dynamical system, we can
assume without loss of generality that the map $T$ is invertible. Now,
fix a time window $[0,L]$ on which we will identify the limit of $X^n$.
By stationarity, we have that
%
\[
X^n(t) = n^{-1/2}\sum_{i=0}^{nt-1}
v \bigl(T^i \bigr) \stackrel{\mathrm {dist}} {=}n^{-1/2}\sum
_{i=0}^{nt-1} v \bigl(T^{i-nL} \bigr).
\]
By setting $i = nL-k$, the above equals
%
\[
\sum_{k=nT-nt}^{nT} v \bigl(T^{-k}
\bigr) = \sum_{k=1}^{nT} v
\bigl(T^{-k} \bigr) - \sum_{k=1}^{n(T-t)}
v \bigl(T^{-k} \bigr).
\]
Now, if we define the \textit{backward time} partial sum $X^n_{-}(t) =
n^{-1/2}\sum_{i=1}^{\lfloor nt \rfloor} v(T^{-i})$ then the above calculation
show that
%
\[
X^n(t) = X^n_{-}(L) - X^n_{-}(L-t).
\]
Under ``sufficient'' mixing conditions on the dynamical system, one can
show that $X^n_-$ is a martingale with respect to the backward time
filtration $\mathcal{F}_t = T^{\lfloor nt \rfloor-1} \mathcal{M}$.
Thus, using the central
limit theorem for ergodic stationary $L^2$ martingale difference
sequences \cite{mcleish74} it follows that $X^n_- \Rightarrow X$ where
$X= D^{1/2}W$. And thus,
%
\[
X^n \Rightarrow X(L) - X(L-t) \stackrel{\mathrm{dist}} {=}X(t),
\]
using the time reversal property of Brownian motion. Now that we know
the limiting behaviour of $X^n$, we can use the tools from \cite
{kurtz91} to identify the limiting behaviour of $\mathbb{X}^n$. In
particular, Theorem~2.2 from \cite{kurtz91} allows us to identify the
limit of integrals against the martingale $X^n_-$, so all we have to do
is rewrite $\mathbb{X}^n$ in backward time, so that it becomes an integral
driven by $dX^n_-$ (plus corrections). Using this idea, we show that
$(X^n,\mathbb{X}^n) \Rightarrow(X,\mathbb{X})$ where
%
\[
\mathbb{X}(t) = \int_0^{t} X(s) \otimes\circ\,
dX(s) + \nu t,
\]
where the integral is of Stratonovich type. All that remains is to
prove the discrete tightness estimates. To do so, we again write the
pair $(X^n,\mathbb{X}^n)$ in terms of the martingale $X^n_-$ and stochastic
integrals driven by $dX^n_-$. Since these are both martingales, we can
apply the BDG inequality and the tightness estimates follow (somewhat)
easily, using the ergodic properties and stationarity of $T$.
It follows from Theorem~\ref{thmmain} that $Y^n \Rightarrow Y$ where
$Y$ solves the SDE
%
\[
dY = h(Y)\circ dW + \nu\dvtx  \mathbb{H}(Y) \,dt.
\]
\upqed\end{pf*}
It is important to note that although the diffusion approximation is
essentially a consequence of the martingale central limit theorem,
those martingale sequences only appear in \textit{backward time}. In
particular, the fast--slow system cannot be written as an equation
driven by a semi-martingale, so the theory of good semi-martingales
cannot be applied directly to the fast--slow system. The advantage of
Theorem~\ref{thmmain} is that even though it is not possible to apply
martingale limit theory to the recursion, it is quite easy to do so to
the \textit{noise} processes $(X^n,\mathbb{X}^n)$.

\begin{rmk}
In the companion paper \cite{kelly14}, the details are far more
complicated than we present above. For example, the backward time
object $X^n_-$ is not in fact a good martingale, but rather a
reasonable perturbation of a good martingale, as described in Example~\ref{eggood}. This makes matters more interesting.
\end{rmk}

\subsection{\texorpdfstring{Connection to \protect\cite{kurtz91}}{Connection to [24]}}

We will now briefly comment on the connection between Theorem~\ref{thmmain} and the tools introduced in \cite{kurtz91}. Given a
partition $\mathcal{P}^n$, consider the recursion
%
\[
Y^n_{k+1} = Y^n_k + V
\bigl(Y^n_k \bigr)\xi^n_k,
\]
where the increments $\xi^n$ are defined in such a way that the
step-function $X^n(t) = \sum_{i=0}^{k-1} \xi^n_i$ were a
semi-martingale with respect to some given sequence of filtrations, the
random walk recursion provides a nice example. It follows that the path
$Y^n(\cdot)$ solves the equation
%
\[
Y^n(t) = Y^n(0) + \int_0^t
V \bigl(Y^n(s-) \bigr)\,dX^n(s),
\]
where the integral is of It\^o type, as defined Section~\ref{srw}.
Suppose moreover that $X^n = M^n + A^n$ where $M^n$ is a good sequence
of semi-martingales and $A^n \Rightarrow0$. Also define
%
\[
H^n(t) = \bigl[M^n,A^n \bigr]_t, \qquad
K^n(t) = \int_0^t
A^n(s) \otimes dA^n(s),
\]
where the integral is of It\^o type. Suppose moreover that
%
\begin{equation}
\label{efours}
\bigl(X^n,A^n,H^n,K^n
\bigr)\Rightarrow (X,0,H,K)
\end{equation}
in the sup-norm topology, where the limits are continuous
semi-martingales. Then \cite{kurtz91}, Theorem~5.1,
states that $Y^n
\Rightarrow Y$ where
%
\begin{equation}
\label{eKPlimit}
dY = V(Y)\,dX + \mathbb{V}\dvtx  d(H-K)
\end{equation}
and the integrals are of It\^o type.

Let us see how this fits into Theorem~\ref{thmmain}. It is not hard to
see that the assumption~\eqref{efours} implies the assumption
$(X^n,\mathbb{X}^n) \Rightarrow(X,\mathbb{X})$. For instance, since
\begin{eqnarray*}
\mathbb{X}^n(t) &=& \int_0^t
X^n(s-) \otimes dX^n(s)
\\
&=& \int_0^t M^n(s-) \otimes
dM^n(s) + \int_0^t
M^n(s-)\otimes dA^n(s)
\\
&&{}+ \int_0^t A^n(s-)\otimes
dM^n(s) + \int_0^t
A^n(s-)\otimes dA^n(s)
\\
&= &\int_0^t M^n(s-) \otimes
dM^n(s) + A^n(t) M^n(t) -
\bigl[M^n,A^n\bigr]_t \\
&&{}+ \int
_0^t A^n(s-)\otimes
dA^n(s)
\\
& = & \int_0^t M^n(s-) \otimes
dM^n(s) + A^n(t) M^n(t) - H^n(t)
+ K^n(t).
\end{eqnarray*}
Combining the fact that $M^n$ is good with \eqref{efours} we see that
$(X^n,\mathbb{X}^n) \Rightarrow(X,\mathbb{X})$ where
\begin{eqnarray*}
\mathbb{X}(t) &= & \int_0^{t} M(s)\otimes dM(s)
-H(t) + K(t)
\\
&=& \int_0^{t} M(s) \otimes\circ\, dM(s) +
\frac{1}{2}[M,M]_t -H(t) + K(t).
\end{eqnarray*}
Taking the tightness estimates for granted, we see that in the case
where $\frac{1}{2}[M,M]_t -H(t) + K(t) = \nu t$, Theorem~\ref{thmmain}
reproduces the diffusion approximation \eqref{eKPlimit}. It is quite
possible to extend the Theorem~\ref{thmmain} so that it only requires, for instance, $[M,M] -H + K$ to be of bounded variation, which would
yield a result closer to that of \cite{kurtz91}, but we do not pursue
this here.

\subsection{Numerical schemes}\label{ssschemes}
Several recent articles have used rough path ideas to study numerical
schemes for stochastic equations. To name a few, \cite
{victoir08,deya12,davie07,riedel14} are all concerned with similar but
typically higher order schemes than \eqref{erecursion}. In \cite
{deya12}, the authors also use the idea that a recursion can be
approximated by an RDE, but only for much higher order Milstein-type
schemes. In a recent preprint \cite{perkowski13}, the authors consider
Euler-type schemes, again by approximating the recursion with a genuine RDE.

The recursion considered in this article handles most numerical schemes
for SDEs, provided the driving noise is a random path with H\"older
exponent $\gamma> 1/3$. It is easy to see that the Euler scheme
%
\[
Y^n_{j+1} = Y^n_{j} + V
\bigl(Y^n_{j} \bigr)X \bigl(j/n,(j+1)/n \bigr)
\]
fits into the framework of \eqref{erecursion}. Using nothing more than
a Taylor expansion, one can also show that another typical numerical
scheme, the semi-implicit scheme, fits into \eqref{erecursion}. This
is defined by
%
\[
Y^n_{j+1} = Y^n_j + \theta V
\bigl(Y^n_j \bigr)X \bigl(j/n,(j+1)/n \bigr) + (1-\theta
)V \bigl(Y^n_{j+1} \bigr)X \bigl(j/n,(j+1)/n \bigr),
\]
where $\theta\in[0,1]$. When $\theta= 1$ this is of course the
(forward) Euler scheme, when $\theta=1/2$ this is the Stratonovich
mid-point scheme and when $\theta=0$ this is the backward Euler scheme.
In \cite{kloden92}, one can find a plethora of schemes that also fit
into the class of recursions defined by \eqref{erecursion}, using
similar arguments to that given below.

In the context of numerical schemes, we see two key areas where the
ability to identify weak limits is beneficial.
\begin{longlist}[2.]
\item[1.] Well-posedness of numerical schemes. When the noise is not a
semi-martingale, it may not be clear whether a limit exists and if it
does---how it should be interpreted. Theorem~\ref{thmmeta} provides a
quick criterion for this situation. In particular, since $X^n=X$ one
need only identify the limit of $\mathbb{X}^n$. If the limit of
$\mathbb{X}^n$
corresponds to a reasonable type of integral (it should correspond to
the method of integration used by the numerical scheme) then the
limiting equation can be interpreted in the sense of that integral.
\item[2.] Numerical schemes that depend on an approximation of the noise,
rather than the exact distribution. Such situations arise if the noise
is difficult to simulate and must instead be approximated, a common
scenario when Gaussianity is not present. One also encounters this
situation in the context of \textit{stochastic climate modeling}, where
ocean--atmosphere equations are driven by an under-resolved source of
noise with persistent correlations in time \cite{penland08} and also in
\textit{data assimilation}, where a perturbation of a stochastic
observation is fed into the numerical simulation of a forecast model.
The article \cite{stuart13} contains a brief overview of the latter idea.
\end{longlist}
Finally, the approximation theory above clearly has applications to
determining the pathwise order of numerical schemes. For example,
suppose that $Y^n$ is defined by the Euler scheme
%
\[
Y^n_{j+1} = Y^n_j + V
\bigl(Y^n_j \bigr)X \bigl(j/n,(j+1)/n \bigr).
\]
Since $X$ does not depend on $n$ the weak limit $Y$ is determined by
the weak limit $\mathbb{X}$ of
%
\[
\mathbb{X}^n(t) = \sum_{j=0}^{\lfloor nt \rfloor-1}
X(0,j/n)\otimes X \bigl(j/n,(j+1)/n \bigr).
\]
Using Theorem~\ref{thmmod} as well as the tools from rough path theory
(Lemma~\ref{thmitomap}) it is easy to show that
%
\[
\sup_{k=0,\ldots,n} \bigl|Y^n \bigl(\tau^n_k
\bigr) - Y \bigl(\tau^n_k \bigr)\bigr| \leq K_{\gamma,n}
\Bigl( \Delta_n^{3\gamma-1} + \sup_{k=0,\ldots,n}\bigl|
\mathbb{X}^n \bigl(\tau ^n_k \bigr) -
\mathbb{X} \bigl(\tau^n_k \bigr)\bigr|^\theta \Bigr),
\]
where $K_{\gamma,n}$ only depends on $n$ through the discrete H\"older
norm of $\mathbb{X}^n$. If $X$ were Brownian motion, then one can trivially
calculate moments of $|\mathbb{X}^n(\tau^n_k) - \mathbb{X}(\tau
^n_k)|$ exactly,
thus obtaining a rate of convergence is simple. However, obtaining the
\textit{optimal} rate of convergence is slightly more subtle. The topic
of convergence rates will not be discussed further in this article but
is the subject of a future article.

\section{A taste of rough path theory}\label{srough}
In this section, we will serve an appetizer in rough path theory. For
the full course, we recommend \cite{frizhairer13}, which is closely
aligned with the exposition below.

\subsection{Space of rough paths}


A rough path has two components to its definition, an algebraic one and
an analytic one. The algebraic component ensures that the objects
$X,\mathbb{X}$ do indeed behave like the increments they hope to
imitate. The
analytic component describes the H\"older condition that is required to
construct solution maps. In the definition below, we \textit{always}
require that the exponent $\gamma> 1/3$. We use the notation
$T^2(\mathbb{R}^d) = \mathbb{R}^d \oplus\mathbb{R}^{d\times d} $
for the step-2
tensor product algebra.

\begin{defn}
We say that $\mathbf{X}\dvtx  [0,T]\times[0,T] \to T^2(\mathbb{R}^d)$ is
a \textit{rough path}
if for $\mathbf{X}= (X,\mathbb{X})$
%
\begin{eqnarray}
 X(s,t) &=& X(s,u) + X(u,t)\quad \mbox{and}
\nonumber
\\[-8pt]
\label{erp2}
\\[-8pt]
\nonumber
\mathbb{X}(s,t) &=& \mathbb{X}(s,u) + \mathbb{X}(u,t) + X(s,u)\otimes X(u,t)
\end{eqnarray}
for all $s,u,t \in[0,T]$. These are known as \textit{Chen's relations}.
If moreover we have that
%
\begin{equation}
\label{erp1}
\bigl|X(s,t)\bigr| \lesssim|t-s|^\gamma\quad \mbox{and}\quad \bigl|
\mathbb{X}(s,t)\bigr| \lesssim|t-s|^{2\gamma}
\end{equation}
for all $s,t \in[0,T]$ then $\mathbf{X}$ is a \textit{$\gamma$-H\"older rough
path}. The set of $\gamma$-H\"older rough paths will be denoted
$\mathcal{C}^\gamma([0,T]; \mathbb{R}^d)$. Every rough path defines
a path $\mathbf{X}\dvtx  [0,T]\to T^2(\mathbb{R}^d)$
by setting $\mathbf{X}(t) = \mathbf{X}(0,t)$. Likewise, we could
equally have defined
rough paths as paths $\mathbf{X}\dvtx  [0,T] \to T^2(\mathbb{R}^d)$ and
then simply taken
Chen's relations as a definition for $\mathbf{X}(s,t)$. This identification
between paths and increments will be used frequently throughout the article.
\end{defn}

We will make use of two metric spaces of rough paths. First,
$\mathcal{C}^\gamma([0,T]; \mathbb{R}^d)$ is a metric space when
furnished with the metric
%
\[
\rho_\gamma(\mathbf{X},\tilde{\mathbf{X}}) = \sup_{s\neq t \in
[0,T]}
\frac
{|X(s,t)-\tilde{X}(s,t)|}{|s-t|^\gamma} + \sup_{s\neq t \in[0,T]} \frac
{|\mathbb{X}(s,t)-\tilde{\mathbb{X}}(s,t)|}{|s-t|^{2\gamma}}.
\]
It is easy to check the interpolation inequality
%
\begin{equation}
\label{einterpolation}
\rho_\alpha\leq\rho_\beta^{\alpha/\beta}
\rho_0^{1-\alpha
/\beta}
\end{equation}
for any $0\leq\alpha\leq\beta$. We also make use of the related
$\gamma$-H\"older ``norm''
%
\[
\vert \!\vert \!\vert\mathbf{X} \vert \!\vert \!\vert_\gamma= \sup
_{s\neq t \in[0,T]} \frac
{|X(s,t)|}{|s-t|^\gamma} + \sup_{s\neq t \in[0,T]}
\frac{|\mathbb{X}
(s,t)|^{1/2}}{|s-t|^\gamma},
\]
which is by definition finite on $\mathcal{C}^\gamma([0,T]; \mathbb
{R}^d)$. Clearly, we have that
%
\[
\bigl|X(s,t)\bigr| \leq\vert \!\vert \!\vert\mathbf{X} \vert \!\vert \!\vert _\gamma|s-t|^\gamma
\quad \mbox{and}\quad \bigl|\mathbb{X}(s,t)\bigr| \leq\vert \!\vert \!\vert\mathbf{X} \vert \!\vert \!\vert_\gamma^2 |s-t|^{2\gamma}
\]
for all $s,t \in[0,T]$ and $\mathbf{X}= (X,\mathbb{X}) \in\mathcal
{C}^\gamma([0,T]; \mathbb{R}^d)$.

The second metric space we make use of is the set of continuous rough
paths $\mathbf{X}\dvtx  [0,T]\to T^2(\mathbb{R}^d)$ endowed with the
uniform metric $\llVert  \cdot \rrVert _\infty$. By this, we simply
mean the sup-norm defined on
functions with range $\mathbb{R}^d \oplus\mathbb{R}^{d\times d}$
(with the
ordinary Euclidean norm on the range). It is easy to see that this
topology is equivalent to that generated by $\rho_0$.
%
\begin{rmk}
There is a good reason for using $\vert \!\vert \!\vert\cdot \vert
\!\vert \! \vert_\gamma$ in addition to
$\rho_\gamma$. This is due to the relationship between the Euclidean
norm on and the Carnot--Caratheodory norm, defined on homogeneous
groups. This will be utilised in Section~\ref{sproperties}.
\end{rmk}

\subsection{Rough differential equations}
%
For $\mathbf{X}\in\mathcal{C}^\gamma([0,T]; \mathbb{R}^d)$ and $V
\in C_b^2$, there is a class of paths
$Y \dvtx  [0,T]\to\mathbb{R}^e$ known as \textit{controlled rough paths}, for
which one can define the integral $\int_0^t V(Y) \,d\mathbf{X}$. We say
that $Y$
is an $X$-controlled rough path if $Y(s,t) = Y(t) - Y(s)$ has the form
%
\[
Y^i(s,t) = Y'_i(s) X(s,t) + O
\bigl(|t-s|^{2\gamma} \bigr)
\]
for all $i = 1 ,\ldots, e$ and $0 \leq s \leq t \leq T$, where $Y'_i
\dvtx
[0,T]\to\mathbb{R}^{e\times d}$ is a $\gamma$-H\"older path. For a
thorough treatment of controlled rough paths and their use in defining
the above integrals, see \cite{frizhairer13}, Section~4.

The integral is defined as a compensated Riemann sum
%
\[
\int_0^t V(Y) \,d\mathbf{X}= \lim
_{\mathcal{P}\to0} S_\mathcal{P},
\]
where
%
\[
S^i_\mathcal{P}= \sum_{[t_k,t_{k+1}] \in\mathcal{P}}
V^i \bigl(Y(t_k) \bigr)X(t_k,t_{k+1})
+ \sum_{j=1}^e \bigl(
Y'_j (t_k) \otimes\partial_j
V^i \bigl(Y(t_k) \bigr) \bigr) \dvtx \mathbb
{X}(t_k,t_{k+1})
\]
and $\mathcal{P}$ denotes a partition of $[0,t]$. Note that the
integral is
defined pathwise, for each $\mathbf{X}\in\mathcal{C}^\gamma([0,T];
\mathbb{R}^d)$.

A controlled rough path $Y$ is said to solve the RDE $dY = V(Y)
\,d\mathbf{X}$
with initial condition $Y(0) = \eta$ if it solves the integral equation
%
\[
Y(t) = \eta+ \int_0^t V(Y) \,d\mathbf{X}
\]
for all $t \in[0,T]$. In this case, we write $Y = \Phi(\mathbf{X})$. When
required, we write $Y = \Phi(\mathbf{X}; V,\eta,s)$ to denote the
solution to
$Y(t) = \eta+ \int_s^t V(Y)\,d\mathbf{X}$ with $t\geq s$. For a thorough
treatment of RDEs, see \cite{frizhairer13}, Section~8.

We will now state a few basis results concerning RDEs, proofs can be
found in \cite{gubinelli04,davie07,frizhairer13}.

\begin{prop}
If $V \in C_b^3$ and $\mathbf{X}\in\mathcal{C}^\gamma([0,T];
\mathbb{R}^d)$ then for each initial
condition $\xi$ and any $s < T $, there exists a unique global solution
$Y = \Phi(\mathbf{X}; V, \xi,s)$. Moreover,
%
\begin{equation}\label{eRDEsoln}
Y(t) = Y(s) + V \bigl(Y(s) \bigr)X(s,t) + \mathbb{V} \bigl(Y(s) \bigr)\dvtx
\mathbb{X}(s,t) + R(s,t),
\end{equation}
for all $s,t \in[0,T]$, where $
|R(s,t)| \lesssim(1\wedge\vert \!\vert \!\vert\mathbf{X} \vert
\!\vert\!  \vert_\gamma^3) |s-t|^{3\gamma} $.
\end{prop}

\begin{rmk}
To see that the remainder scales in this particular way, see the proof
of \cite{davie07}, Lemma~3.4.
\end{rmk}

%
%

\begin{lemma}\label{lemic}
If $Y = \Phi(\mathbf{X};V,\eta)$ and $\tilde{Y}= \Phi(\mathbf
{X};V,\tilde{\eta})$, then
$|Y(t) - \tilde{Y}(t)|\lesssim\vert \!\vert \!\vert\mathbf{X} \vert\!
  \vert\!  \vert_\gamma|Y(s) - \tilde{Y}(s)|$
for any $s\leq t\leq T$, where the implied constant depends only on $T,V$.
\end{lemma}


\begin{lemma}\label{thmitomap}
Suppose that $V \in C^3_b$ and $\mathbf{X},\tilde{\mathbf{X}}\in
\mathcal{C}^\gamma([0,T]; \mathbb{R}^d)$
satisfying $\rho_\gamma(\mathbf{X},0), \rho_\gamma(\tilde{\mathbf
{X}},0) \leq M$. Then,
on any time window $[0,T]$, the map $\Phi(\cdot)$ satisfies the
following local Lipschitz estimate
%
\[
\bigl\llVert \Phi(\mathbf{X})-\Phi(\tilde{\mathbf{X}}) \bigr\rrVert
_{\infty} \leq C_M\rho_\gamma(\mathbf{X},\tilde{
\mathbf{X}}),
\]
where $C_M$ depends only on $M$ and $T$.
\end{lemma}

The next lemma is a slight modification of a result in \cite{friz10},
hence we only sketch the proof.

\begin{lemma}
Let $0 \leq\gamma< \alpha$. Then the ball
%
\[
B_{R,\alpha} = \bigl\{ \mathbf{X}\in\mathcal{C}^\alpha \bigl([0,T];
\mathbb {R}^d \bigr)\dvtx \vert \!\vert \!\vert\mathbf{X} \vert \!\vert \!\vert
_{\alpha} \leq R \bigr\}
\]
is compact in the space $\mathcal{C}^\gamma([0,T]; \mathbb{R}^d)$.
\end{lemma}
\begin{pf}
The proof is a standard modification of a similar statement found in
\cite{friz10}. Fix a sequence $\{\mathbf{X}^n\} \subset B_{R,\alpha
}$. Use
Arzela--Ascoli to find a subsequence that converges\vspace*{1pt} in the sup-norm
topology. Use the interpolation \eqref{einterpolation} between
$\mathcal{C}^\alpha\subset\mathcal{C}^\gamma\subset\mathcal{C}^0 $
to show that this subsequence also converges in $\mathcal{C}^\gamma$.
Since $\mathcal{C}^\gamma$ is a metric space, sequential compactness
implies compactness.
\end{pf}

The final result, which is a direct corollary of \cite{friz10}, Theorem~17.3, allows us to translate RDE solutions to Stratonovich SDEs.

\begin{lemma}\label{lemrpsde}
Suppose that $V \in C^3_b$, $\mathbf{X}= (X,\mathbb{X}) \in\mathcal
{C}^\gamma([0,T]; \mathbb{R}^d)$ and let $Y
= \Phi(\mathbf{X}; V)$. Suppose that $X$ is a continuous
semi-martingale and
that $\mathbb{X}$ can be written
%
\begin{equation}
\label{eitoplus}
\mathbb{X}^{\alpha\beta}(t) = \int_0^t
X^\beta(s) \circ dX^\alpha (s) + \nu ^{\alpha\beta} t,
\end{equation}
for $\alpha, \beta= 1 ,\ldots, d$. where the integral\vspace*{1pt} on the right-hand
side is defined in the Stratonovich sense and where $\nu\in\mathbb{R}
^{d\times d}$. Then $Y$ satisfies the SDE
%
\begin{equation}
\label{eito}
dY = V(Y)\circ dX + \mathbb{V}(Y) \dvtx \nu\, dt.
\end{equation}
\end{lemma}

\begin{rmk}
In all of the results of this section, if in addition to $\gamma> 1/3$,
it is also known that $\gamma$ can be taken arbitrarily close to $1/2$
(as for Brownian rough paths), then the condition $V\in C^{3} $ can be
relaxed to $V \in C^{2+}$. For instance, the versions of the results in
any of \cite{gubinelli04,davie07,friz10} will adhere to this.
\end{rmk}

\section{Rough path recursions}\label{sdrp}
In this section, we introduce \textit{rough path recursions}, driven by
\textit{rough step functions}. Before proceeding with the definitions, we
must introduce some assumptions and terminology.
\subsection{Partitions of $[0,T]$}
Fix an interval $[0,T]$ and let $\mathcal{P}_n = \{\tau^n_k \dvtx  k=0
,\ldots, N_n \}$
be a partition of $[0,T]$, that is $0 = \tau_0^n \leq\tau_1^n \leq
\cdots\leq\tau_{N_n}^n = T $. We also introduce the \textit{mesh size}
$\Delta_n = \max_{k} |\tau_{k+1}^n - \tau_k^n| $. For all the
partitions considered in this article, we will assume that
%
\begin{equation}
\label{eassp}
\Delta_n \to0\qquad \mbox{as $n\to\infty$} \quad\mbox{and}\quad \sup
_{n \geq1} N_n \Delta_n < \infty.
\end{equation}
The first assumption is obviously natural, the second condition is
effectively saying that the largest bin the  partition does not
shrink too much slower than the smallest bin in the partition. Given
some $u \in[0,T]$, we will also use the notation $\tau^n(u)$ to denote
the largest mesh point $\tau^n(u) \in\mathcal{P}_n$ with $\tau^n(u)
\leq u$.
It follows from \eqref{eassp} that $\tau^n(u) \to u$ as $n\to\infty$.

\subsection{Rough path recursions}
We will now define rough step functions and rough path recursions rigorously.

\begin{defn}[(Rough step functions)]
Fix a partition $\mathcal{P}_n$ of $[0,T]$ and suppose $\xi^n_j \in
\mathbb{R}^d$
and $\Xi^n_j \in\mathbb{R}^{d\times d}$ for all $j=0, \ldots, N_n-1$. The
\textit{rough step function} above the increments $(\xi^n,\Xi^n)$ is a
path $\mathbf{X}^n = (X^n,\mathbb{X}^n) \dvtx  [0,T]\to T^2(\mathbb
{R}^d)$ defined by
%
\[
X^n(t) = \sum_{j=0}^{k-1}
\xi^n_{j} \quad\mbox{and}\quad \mathbb {X}^n(t) = \sum
_{i=0}^{k-1} \sum
_{j=0}^{i-1} \xi^n_i \otimes
\xi^n_j + \sum_{i=0}^{k-1}
\Xi^n_i,
\]
where $\tau^n_k = \tau^n(t)$. We similarly define the incremental paths
%
\[
X^n(s,t) = \sum_{j=l}^{k-1}
\xi^n_{j} \quad\mbox{and} \quad\mathbb {X}^n(s,t) = \sum
_{i=l}^{k-1} \sum
_{j=l}^{i-1} \xi^n_i \otimes
\xi^n_j + \sum_{i=l}^{k-1}
\Xi^n_i,
\]
where $\tau^n_k = \tau^n(t)$ and $\tau^n_l = \tau^n(s)$. We will often
employ the shorthand $X^n(\tau^n_j, \tau^n_k) = X^n_{j,k}$ and
$\mathbb{X}
^n(\tau^n_j,\tau^n_k) = \mathbb{X}^n_{j,k}$. We define the discrete
$\gamma
$-H\"older ``norm''  $\vert \!\vert \!\vert\cdot \vert \! \vert\!
\vert_{\gamma,n}$ by
%
\[
\bigl\vert \!\bigl\vert \!\bigl\vert\mathbf{X}^n \bigr\vert \!\bigr\vert \!\bigr\vert_{\gamma
,n}
\stackrel{{\mathrm{def}}} {=}\max_{\tau^n_j,\tau^n_k \in
\mathcal{P}_n} \frac
{|X^n_{j,k}|}{|\tau^n_j-\tau^n_k|^\gamma} +
\max_{\tau^n_j,\tau
^n_k \in
\mathcal{P}_n} \frac{|\mathbb{X}^n_{j,k}|^{1/2}}{|\tau^n_j-\tau
^n_k|^\gamma}.
\]
In particular, we see that
%
\[
\bigl|X^n \bigl(\tau^n_j,\tau^n_k
\bigr)\bigr| \leq\bigl\vert \!\bigl\vert \!\bigl\vert\mathbf{X}^n \bigr\vert \!\bigr\vert \!\bigr\vert_{\gamma,n}\bigl|\tau^n_j - \tau
^n_k\bigr|^\gamma\quad \mbox{and}\quad \bigl|\mathbb{X}^n
\bigl(\tau^n_j,\tau ^n_k
\bigr)\bigr| \leq \bigl\vert \!\bigl\vert \!\bigl\vert\mathbf{X}^n \bigr\vert \!\bigr\vert \!\bigr\vert_{\gamma
,n}^2\bigl| \tau^n_j -
\tau^n_k\bigr|^{2\gamma}
\]
for all {mesh points} $\tau^n_j,\tau^n_k \in\mathcal{P}_n$. Since
it is a
maximum over a finite set, the discrete H\"older norm is finite for
every fixed $n$. It will only play a role in an asymptotic sense.
\end{defn}

\begin{defn}[(Rough path recursions)]\label{drpr}
Fix a sequence of partitions $\{\mathcal{P}_n\}_{n\geq1}$. A \textit{rough path
recursion} $\{Y^n\}_{n\geq1}$ with $Y^n \dvtx  [0,T]\to\mathbb{R}^e$ is defined
by $Y^n(t) = Y^n_j$, where $\tau^n_j = \tau^n(t)$ and satisfying the recursion
%
\begin{equation}
\label{erps}
Y^n_{ j+1} = Y^n_{j}
+ V \bigl(Y^n_{ j} \bigr)\xi^{n}_{j}
+ \mathbb{V} \bigl(Y^n_{j} \bigr) \dvtx  \Xi
^{n}_{j}+r^n_{j},
\end{equation}
for all $j=0,\ldots, N_n-1$ with $Y^n_0 = \eta$ and arbitrary $\xi^n_j
\in\mathbb{R}^d$ and $\Xi^n_j \in\mathbb{R}^{d\times d}$. The
remainder is
assumed to satisfy the estimate
%
\begin{equation}
\label{eRPRremainder}
|r_j| \lesssim \bigl(1\wedge\bigl\vert \!\bigl\vert \!\bigl\vert
\mathbf{X}^n \bigr\vert \!\bigr\vert \!\bigr\vert_{\gamma,n}^3 \bigr)
\Delta_n^{3\gamma},
\end{equation}
where the implied constant is uniform in $n$ and $\mathbf{X}^n$ is the rough
step function over the increments $(\xi^n,\Xi^n)$. We will use the
notation $Y^n = \Phi^n(\mathbf{X}^n)$.
\end{defn}

\begin{rmk}
The estimate \eqref{eRPRremainder} is picked as it seems to be the
most naturally occurring upper bound in applications. However, at the
expense of a few extra constraints we could equally use
$|r^n_j|\lesssim(1\wedge\phi(C_n))\Delta^\lambda_n$ for any
increasing function $\phi\dvtx  [0,\infty) \to[0,\infty)$ and $\lambda>1$.
\end{rmk}

\begin{rmk}\label{rmkinterpolation}
Although we require that $Y^n$ be constant in between mesh points, it
is easy to see that all the properties of rough path recursions
discussed in the sequel are still true if we assume that $Y^n$ is
defined by a reasonable interpolation between mesh points. For example,
even though the solution to an RDE satisfies the recursion \eqref
{erps}, it is not a rough path recursion. However, it can be well
approximated by a rough path recursion, and any convergence results for
rough path recursion easily imply convergence results for the
associated RDE solution. As a more general course of action, we could
have defined a rough path recursion to any path satisfying \eqref
{erps} as well as
%
\[
\bigl|Y^n(t) - Y^n \bigl(\tau^n(t) \bigr)\bigr|
\lesssim D_n \Delta_n^\mu,
\]
for some sequence of constants $D_n$ and $\mu>0$. It is clear that,
assuming the right conditions on $D_n$ and $\mu$, all the statements
made in the sequel regarding rough path recursions are unaltered if we
were to adopt this more general definition.
\end{rmk}

\section{Properties of rough path recursions}\label{sproperties}
In this section, we discuss some useful properties of rough step
functions and their associated rough path recursions. The main result,
Lemma~\ref{thmytilde}, states that every rough path recursion can be
approximated arbitrarily well by the solution to a rough differential
equation. At the heart of this result is the fact that for every rough
step function $\mathbf{X}^n$, one can find a genuine rough path
$\tilde{\mathbf{X}}^n$
that agrees with $\mathbf{X}^n$ on $\mathcal{P}_n$. This is the
content of Lemma~\ref{thmbxtilde}. Before stating the theorems, we must introduce some
terminology associated with \textit{geometric} rough paths. For a more
detailed exposition of this material, see \cite{friz10} and also
\cite{frizhairer13}, Section~2.2.

We first define $G^2(\mathbb{R}^d)$, the step-2 nilpotent Lie group, by
$G^2(\mathbb{R}^d) =  \exp (\mathcal{G}^2(\mathbb{R}^d)  )$
where $\mathcal{G}^2(\mathbb{R}^d)$ is the step-2 Lie algebra over
$\mathbb{R}^d$ and
where $\exp$ is the tensor exponential. In particular, $g \in
G^2(\mathbb{R}
^d)$ if and only if
%
\[
g = \exp \bigl( a_\alpha e_\alpha+ b_{\beta,\kappa}
[e_\beta,e_\kappa] \bigr) = a_\alpha e_\alpha+
\bigl(\tfrac{1}{2}a_\beta a_\kappa+ b_{\beta,\kappa} -
b_{\kappa,\beta} \bigr)e_\beta\otimes e_\kappa,
\]
where $\{e_\alpha\}$ denotes the canonical basis of $\mathbb{R}^d$
and we
employ the Einstein summation convention. It is easy to see that every
element $A \in T^2(\mathbb{R}^d)$ can be decomposed into
%
\[
A = g + z,
\]
where $g \in G^2(\mathbb{R}^d)$ and $z \in \operatorname{Sym}(\mathbb{R}^{d
\times
d})$. The pair $(G^2(\mathbb{R}^d), \otimes)$ forms a group. This
group has
a \textit{homogeneous metric} known as the Carnot--Caratheodory metric,
defined using geodesic paths. To make this precise, we first define the
\textit{signature} of a smooth path. Let $\mathit{BV}([0,1];\mathbb{R}^d)$ be the
space of paths $\Gamma\dvtx  [0,1]\to\mathbb{R}^d$ with bounded
variation. For
$\Gamma\in \mathit{BV} ([0,1];\mathbb{R}^d)$, the signature is defined by
%
\[
\mathbf{S}(\Gamma) (s,t) = \biggl(\Gamma(t)-\Gamma(s), \int
_s^t \bigl(\Gamma(u) - \Gamma(s) \bigr)\otimes
d \Gamma(u) \biggr),
\]
where the integral is constructed in the Riemann--Stieltjes sense. The
Carnot--Caratheodory (CC) norm is defined by
%
\[
\llVert g \rrVert _{\mathrm{CC}} \stackrel{\mathrm{def}} {=}\inf \biggl\{
\int_0^1 |\,d\Gamma| \dvtx  \Gamma\in \mathit{BV} \bigl([0,1];
\mathbb{R}^d \bigr) \mbox{ and } g = \mathbf{S}(\Gamma) (0,1) \biggr\}.
\]
The following result (which is a refinement of Chow's theorem) shows
that the norm is well defined (see \cite{friz10}, Theorem~7.32, for a
simple proof).

\begin{lemma}\label{lemgeodesic}
If $g \in G^2(\mathbb{R}^d)$ then there exists a path $\Gamma\dvtx  [u,v]
\to
\mathbb{R}^d$ with $|\dot{\Gamma}|=\mbox{const}$ such that $g = \mathbf
{S}(\Gamma)(u,v)$ and
$\llVert  g \rrVert _{\mathrm{CC}} = (v-u)|\dot{\Gamma}|$.
\end{lemma}

The CC norm can also be ``compared'' with the usual Euclidean norms in
the following way. Suppose that $g \in G^2(\mathbb{R}^d)$ can be decomposed
into $g = g_1 + g_2$ where $g_1 \in\mathbb{R}^d$ and $g_2 \in\mathbb{R}
^{d\times d}$, then we have the comparison
%
\begin{equation}
\label{eCCeq} |g_1| + |g_2|^{1/2} \lesssim
\llVert g \rrVert _{\mathrm{CC}} \lesssim|g_1| +
|g_2|^{1/2}.
\end{equation}
This comparison will be useful in the sequel.
%
%

\begin{lemma}\label{thmbxtilde}
Let $\{\mathbf{X}^n\}_{n\geq1}$ be a rough step function on a
partition $\{\mathcal{P}
_n\}_{n\geq1}$ and let $\gamma\in(1/3,1/2]$. For each $n$, there
exists $\tilde{\mathbf{X}}^n \in\mathcal{C}^\gamma([0,T]; \mathbb
{R}^d)$ with
%
\[
\tilde{\mathbf{X}}^n \bigl(\tau^n_j \bigr)
= \mathbf{X}^n \bigl(\tau^n_j \bigr)
\qquad\mbox{for all $\tau^n_j\in \mathcal{P}_n$}.
\]
Moreover, we have that $\vert \!\vert \!\vert\tilde{\mathbf{X}}^n
\vert  \!\vert \! \vert_\gamma\lesssim\vert \!\vert \!\vert\mathbf
{X} ^n \vert \! \vert\!  \vert_{\gamma,n}$ where the implied constant
is uniform in $n$.
\end{lemma}
\begin{pf}
We will start by constructing $\tilde{\mathbf{X}}^n(s,t)$ for fixed
$s,t\in[\tau
^n_j,\tau^n_{j+1}]$ with $s\leq t$. First, note that since $\mathbf{X}
^n_{j,j+1} \in T^2(\mathbb{R}^d)$, we have the decomposition
%
\[
\mathbf{X}^n_{j,j+1}= g \bigl(\tau^n_j,
\tau^n_{j+1} \bigr) + z \bigl(\tau^n_j,
\tau ^n_{j+1} \bigr),
\]
where $g(\tau^n_j,\tau^n_{j+1}) \in G^2(\mathbb{R}^d)$ and $ z(\tau
^n_j,\tau
^n_{j+1}) \in\operatorname{Sym}(\mathbb{R}^{d\times d})$ is defined by
%
\[
z \bigl(\tau^n_j,\tau^n_{j+1}
\bigr) = \tfrac{1}{2} \bigl( \mathbb{X}^{n,\alpha
\beta
}_{j,j+1} +
\mathbb{X}^{n,\beta\alpha}_{j,j+1} - X^{n,\alpha
}_{j,j+1}X^{n,\beta}_{j,j+1}
\bigr) e_\alpha\otimes e_\beta.
\]
We will now define
%
\[
\tilde{\mathbf{X}}^n(s,t)= g(s,t) + z(s,t).
\]
First, define $z(s,t)$ by a simple linear interpolation
%
\[
z(s,t) = \frac{t-s}{\tau^n_{j+1}-\tau^n_j}z \bigl(\tau^n_j,\tau
^n_{j+1} \bigr).
\]
Now, to define $g(s,t)$, we know from Lemma~\ref{lemgeodesic} that
there exists a path $\Gamma\dvtx [\tau^n_j,\tau^n_{j+1}] \to\mathbb
{R}^d$ such
that $|\dot{\Gamma}|=\mbox{const}$ and $g(\tau^n_j,\tau^n_{j+1}) =
S(\Gamma)(\tau
^n_j,\tau^n_{j+1})$ and $\llVert  g(\tau^n_j,\tau^n_{j+1}) \rrVert
_{\mathrm{CC}} = (\tau
^n_{j+1}-\tau^n_j)|\dot{\Gamma}|$. We set $g(s,t) = \mathbf
{S}(\Gamma)(s,t)$.

We will now define $\tilde{\mathbf{X}}^n(s,t)$ for arbitrary $s,t \in
[0,T]$ with
$s\leq t$. Suppose without loss of generality that $s\leq\tau^n_j
\leq
\tau^n_k \leq t$ where $\tau^n_{j-1} = \tau^n(s)$ and $\tau^n_k =
\tau
^n(t)$. Then we define $ \tilde{\mathbf{X}}^n(s,t) = (\tilde
{X}^n,\tilde{\mathbb{X}}
^n)(s,t)$ using Chen's relations
%
\[
\tilde{X}^n(s,t) = \tilde{X}^n \bigl(s,
\tau^n_j \bigr)+ X^n_{j,k} + \tilde
{X}^n \bigl(\tau^n_k,t \bigr)
\]
and
\begin{eqnarray*}
\tilde{\mathbb{X}}^n(s,t) &= &\tilde{\mathbb{X}}^n
\bigl(s,\tau^n_j\bigr) + \mathbb{X}^n_{j,k}
+ \tilde{\mathbb{X}} ^n\bigl(\tau^n_k,t
\bigr)
\\
&&{}+ X^n_{j,k}\otimes\tilde{X}^n\bigl(
\tau^n_k,t\bigr) + \tilde{X} ^n\bigl(s,
\tau^n_j\bigr)\otimes\tilde{X}^n\bigl(
\tau^n_j,t\bigr).
\end{eqnarray*}
We will now check that $\tilde{\mathbf{X}}^n$ satisfies the
requirements of the
theorem. First, we will show that Chen's relations hold. It is easy to
see that Chen's relations hold when restricted to the interval $s,u,t
\in[\tau^n_j,\tau^n_{j+1}]$. Indeed, since $g(s,t)$ is a signature and
since $z(s,t)$ is an increment, both objects individually obey Chen's
relation. Thus, if we write $g = g_1 + g_2$, with $g_1 \in\mathbb{R}^d$
and $g_2 \in\mathbb{R}^d\otimes\mathbb{R}^d$ then
%
\[
\tilde{\mathbf{X}}^n = \bigl(\tilde{X}^n,\tilde{
\mathbb{X}}^n \bigr) = (g_1, g_2 + z).
\]
Therefore,
%
\[
\tilde{X}^n(s,t) = g_1(s,t) = g_1(s,u) +
g_1(u,t) = \tilde{X}^n(s,u) + \tilde{X}^n(u,t)
\]
and
\begin{eqnarray*}
\mathbb{X}^n(s,t) &=& g_2(s,t) + z(s,t)
\\
&=& g_2(s,u) + g_2(u,t) + g_1(s,u)\otimes
g_1(u,t) + z(s,u) + z(u,t)
\\
&=& (g_2+z) (s,u) + (g_2+z) (u,t) + g_1(s,u)
\otimes g_1(u,t)
\\
&=& \tilde {\mathbb{X}} ^n(s,u) + \tilde{\mathbb{X}}^n(u,t)
+ \tilde{X}^n(s,u)\otimes\tilde {X}^n(u,t),
\end{eqnarray*}
as required. Now, for arbitrary $s\leq u \leq t$ it follows immediately
from the construction that $(\tilde{X}^n,\tilde{\mathbb{X}}^n)$
satisfies Chen's relations.

Using the shorthand $C_n = \vert \!\vert \!\vert\mathbf{X}^n \vert\!
\vert \! \vert_{\gamma,n}$, we will now prove
that
$\vert \!\vert \!\vert\tilde{\mathbf{X}}^n \vert \! \vert\!
\vert_\gamma\lesssim C_n$. First, suppose that
$s,t \in[\tau^n_j,\tau^n_{j+1}]$ with $s\leq t$. Then using the
comparison \eqref{eCCeq} and the construction of $g$, we have\vspace*{-2pt} that
\begin{eqnarray*}
\bigl|\tilde{X}^n(s,t)\bigr| &=& \bigl|g_1(s,t)\bigr| \leq\bigl|g_1(s,t)\bigr|
+ \bigl|g_2(s,t)\bigr|^{1/2}
\\[-2pt]
&\lesssim & \bigl\llVert g(s,t) \bigr\rrVert _{\mathrm{CC}} = \frac{(t-s)}{(\tau
^n_{j+1}-\tau^n_j)}
\bigl\llVert g\bigl(\tau^n_j,\tau^n_{j+1}
\bigr) \bigr\rrVert _{\mathrm{CC}}
\end{eqnarray*}
and again by \eqref{eCCeq} we have\vspace*{-2pt} that
%
\begin{eqnarray}
 \bigl\llVert g\bigl(\tau^n_j,
\tau^n_{j+1}\bigr) \bigr\rrVert _{\mathrm{CC}} &
\lesssim & \bigl|g_1\bigl(\tau ^n_j,\tau
^n_{j+1}\bigr)\bigr| + \bigl|g_2\bigl(
\tau^n_j,\tau^n_{j+1}
\bigr)\bigr|^{1/2}\nonumber\\[-2pt]
\label{efeb9a}
&\lesssim & \bigl|X^n_{j,j+1}\bigr| + \bigl(\bigl|\mathbb{X}^n_{j,j+1}\bigr|
+ \bigl|X^n_{j,j+1}\bigr|^2 \bigr)^{1/2}
\\[-2pt]
& \lesssim &  C_n \bigl(\tau^n_{j+1} -
\tau^n_j \bigr)^\gamma.\nonumber
\end{eqnarray}
It follows that
%
\[
\bigl|\tilde{X}^n(s,t)\bigr| \lesssim C_n \frac{(t-s)}{(\tau^n_{j+1}-\tau
^n_j)^{1-\gamma}} =
C_n (t-s)^\gamma \biggl(\frac{t-s}{\tau
^n_{j+1}-\tau
^n_j}
\biggr)^{1-\gamma} \leq C_n (t-s)^\gamma,
\]
where in the last inequality we use the fact that $\frac{t-s}{\tau
^n_{j+1}-\tau^n_j} \leq1$. By a similar argument, we can show that
%
\begin{equation}
\label{efeb10a}
\quad\bigl|\tilde{\mathbb{X}}^n(s,t)\bigr| \lesssim \biggl(
\frac{t-s}{\tau
^n_{j+1}-\tau
^n_j} \biggr)^2 \bigl( \bigl|X^n_{j,j+1}\bigr|^2
+ \bigl|\mathbb{X}^n_{j,j+1}\bigr| \bigr) + \biggl(\frac{t-s}{\tau^n_{j+1}-\tau^n_j}
\biggr) \bigl|\mathbb{X}^n_{j,j+1}\bigr|\hspace*{-20pt}
\end{equation}
and hence
%
\[
\bigl|\tilde{\mathbb{X}}^n(s,t)\bigr| \lesssim C_n^2
(t-s)^{2\gamma}.
\]
Now suppose $s,t \in[0,T]$ with $s\leq\tau^n_j \leq\tau^n_k \leq t$
as above. By Chen's relations, we have that
\begin{eqnarray*}
\bigl|\tilde{X}^n(s,t)\bigr| &\leq& \bigl|\tilde{X}^n\bigl(s,
\tau^n_j\bigr)\bigr|+ \bigl|X^n_{j,k}\bigr| + \bigl|
\tilde{X} ^n\bigl(\tau^n_k,t\bigr)\bigr|
\\
&\lesssim  & C_n \bigl(\bigl|\tau^n_j-s\bigr|^\gamma+
\bigl|\tau^n_j - \tau ^n_k\bigr|^\gamma+
\bigl|\tau^n_k - t\bigr|^\gamma\bigr) \leq
C_n|t-s|^\gamma
\end{eqnarray*}
and
\begin{eqnarray*}
\bigl|\tilde{\mathbb{X}}^n(s,t)\bigr| &\leq & \bigl|\tilde{\mathbb{X}}^n
\bigl(s,\tau ^n_j\bigr)\bigr| + \bigl|\mathbb{X}^n_{j,k}\bigr|
+ \bigl|\tilde{\mathbb{X}}^n\bigl(\tau^n_k,t
\bigr)\bigr|
\\
&&{}+ \bigl|X^n_{j,k}\otimes\tilde {X}^n\bigl(
\tau^n_k,t\bigr)\bigr| + \bigl|\tilde{X}^n\bigl(s,
\tau^n_j\bigr)\otimes\tilde{X}^n\bigl(
\tau^n_j,t\bigr)\bigr|
\\
&\lesssim & C_n^2 \bigl( \bigl|\tau^n_j-s\bigr|^{2\gamma}
+ \bigl|\tau^n_j - \tau ^n_k\bigr|^{2\gamma}
+ \bigl|\tau^n_k - t\bigr|^{2\gamma}
\\
&&\hspace*{15pt}{}+ \bigl|\tau^n_j-\tau^n_k\bigr|^\gamma\bigl|
\tau^n_k - t\bigr|^\gamma+ \bigl|\tau^n_j-s\bigr|^\gamma
\bigl|\tau^n_j-t\bigr|^\gamma \bigr) \lesssim
C_n |t-s|^{2\gamma}.
\end{eqnarray*}
This completes the proof.
\end{pf}
We now have all the tools we need to prove the main result of this
section, namely that rough path recursions can be well approximated by
the solution to an RDE.

\begin{lemma}\label{thmytilde}
Let $Y^n = \Phi^n(\mathbf{X}^n)$. Let $\tilde{\mathbf{X}}^n$ be any
rough path
satisfying the conditions of Lemma~\ref{thmbxtilde} and let $\tilde{Y}^n
= \Phi(\tilde{\mathbf{X}}^n)$. Then
%
\[
\bigl\llVert \tilde{Y}^n - Y^n \bigr\rrVert
_{\infty} \lesssim \bigl(1\wedge\bigl\vert \!\bigl\vert \!\bigl\vert\mathbf{X}^n
\bigr\vert \!\bigr\vert \!\bigr\vert_{\gamma
,n}^4 \bigr) \Delta_n^{3\gamma-1},
\]
for any $\gamma\in(1/3,1/2]$, where the implied constant is uniform
in $n$.
\end{lemma}

\begin{pf}
We will again use the shorthand $C_n = \vert \!\vert \!\vert\mathbf
{X}^n \vert \! \vert\!  \vert_{\gamma,n}$. For
any $t \in[0,T]$, we have that
%
\[
\bigl|\tilde{Y}^n(t) - Y^n(t) \bigr| \leq\bigl|\tilde{Y}^n(t)
- \tilde{Y}^n \bigl(\tau ^n_k \bigr) \bigr| + \bigl|
\tilde{Y}^n \bigl(\tau^n_k \bigr) -
Y^n \bigl(\tau^n_k \bigr)\bigr|,
\]
where $\tau^n_k = \tau^n(t)$, and hence $Y^n(t) = Y^n(\tau^n_k)$. It
follows from \eqref{eRDEsoln} that
%
\[
\tilde{Y}^n(t) = \tilde{Y}^n \bigl(\tau^n_k
\bigr) + V \bigl(\tilde{Y}^n \bigl(\tau^n_k
\bigr) \bigr) \tilde{X} ^{n} \bigl(\tau^n_k,t
\bigr) + \mathbb{V} \bigl(\tilde{Y}^n \bigl(\tau^n_k
\bigr) \bigr) \dvtx \tilde {\mathbb{X}}^{n} \bigl(\tau
^n_k,t \bigr) + R \bigl(\tau^n_k,t
\bigr),
\]
where $|R(\tau^n_k,t)|\lesssim\vert \!\vert \!\vert\tilde{\mathbf
{X}}^n \vert \! \vert\!  \vert^3_\gamma(t -\tau
^n_j)^{3\gamma} \leq C_n^3 \Delta_n^{3\gamma}$. Since $\gamma\leq
1/2$, it follows that
%
\[
\bigl|\tilde{Y}^n(t) - \tilde{Y}^n \bigl(
\tau^n_k \bigr)\bigr| \lesssim \bigl(1\wedge
C_n^3 \bigr) \Delta _n^{\gamma}
\lesssim \bigl(1\wedge C_n^4 \bigr) \Delta_n^{3\gamma-1}.
\]
To estimate $|\tilde{Y}^n(\tau^n_k)-Y^n(\tau^n_k)|$ we need some new
terminology. For each $l\leq k$, define $Z^{(l)}_k$ by
%
\[
Z^{(l)}_k \stackrel{{\mathrm{def}}} {=}\Phi \bigl(\tilde{
\mathbf{X}}^n ;V, Y^n_l,
\tau^n_l \bigr) \bigl(\tau^n_k
\bigr).
\]
That is, $Z^{(l)}_k = \hat{Y}(\tau^n_k)$, where $\hat{Y}$ is the unique
solution to the RDE driven by $\tilde{\mathbf{X}}^n$ with vector
field $V$,
initialised at time $\tau^n_l$ with initial condition $Y^n_l$
[recalling that $Y^n_l = Y^n(\tau^n_l)$, as in Definition~\ref{drpr}].
In particular, we have that $\tilde{Y}^n(\tau^n_k) = Z^{(0)}_k$,
$Y^n_k =
Z^{(k)}_k$ and
%
\begin{equation}
\label{eZstep}
\qquad Z^{(k)}_{k+1} = Y^n_k
+ V \bigl(Y^n_k \bigr)\tilde{X}^{n} \bigl(
\tau^n_k,\tau^n_{k+1} \bigr) +
\mathbb{V} \bigl(Y^n_k \bigr)\tilde{\mathbb{X}}^{n}
\bigl(\tau^n_k,\tau^n_{k+1} \bigr)
+ R \bigl(\tau^n_k,\tau^n_{k+1}
\bigr)
\end{equation}
for any $k$. It follows that
%
\begin{equation}
\label{eZ1}
\bigl|\tilde{Y}^n\bigl(\tau^n_k
\bigr) - Y^n\bigl(\tau^n_k\bigr)\bigr| =
\bigl|Z^{(0)}_k - Z^{(k)}_k\bigr| \leq \sum
_{l=0}^{k-1} \bigl|Z^{(l)}_k
- Z^{(l+1)}_k\bigr|.
\end{equation}
But since
%
\[
Z^{(l)}_k = \Phi \bigl(\tilde{\mathbf{X}}^n;
V,Y^n_l,\tau^n_l \bigr) \bigl(
\tau^n_k \bigr) = \Phi \bigl(\tilde{\mathbf{X}}
^n; V, Z^{(l)}_{l+1},\tau^n_{l+1}
\bigr) \bigl(\tau^n_k \bigr)
\]
and $Z^{(l+1)}_k = \Phi(\tilde{\mathbf{X}}^n; V,Z^{(l+1)}_{l+1},\tau
^n_{l+1})(\tau
^n_k)$, it follows from Lemma~\ref{lemic} that
%
\begin{equation}
\label{eZ2}
\qquad\bigl|Z^{(l)}_k - Z^{(l+1)}_k\bigr|
\lesssim \bigl(1\wedge\bigl\vert \!\bigl\vert \!\bigl\vert \tilde{\mathbf{X}}^n \bigr\vert \!\bigr\vert \!\bigr\vert_\gamma \bigr) \bigl|Z^{(l)}_{l+1} -
Z^{(l+1)}_{l+1}\bigr| \leq(1\wedge C_n)
\bigl|Z^{(l)}_{l+1} - Z^{(l+1)}_{l+1}\bigr|.\hspace*{-11pt}
\end{equation}
By \eqref{eZstep}, we have
\begin{eqnarray*}
&& Z^{(l)}_{l+1} - Z^{(l+1)}_{l+1}
\\
&&\qquad=  Y^n_l + V\bigl(Y^n_l\bigr)
\tilde {X}^{n}\bigl(\tau ^n_l,
\tau^n_{l+1}\bigr) + \mathbb{V}\bigl(Y^n_l
\bigr)\dvtx  \tilde{\mathbb{X}}^{n}\bigl(\tau ^n_l,
\tau^n_{l+1}\bigr) + R\bigl(\tau^n_l,
\tau^n_{l+1}\bigr) - Y^n_{l+1}
\\
&&\qquad= Y^n_l + V\bigl(Y^n_l
\bigr)X^{n}_{l,l+1} + \mathbb{V}\bigl(Y^n_l
\bigr)\dvtx  \mathbb {X}^{n}_{l,l+1} + R\bigl(\tau^n_l,
\tau^n_{l+1}\bigr) - Y^n_{l+1}\\
&& \qquad = R
\bigl(\tau^n_l,\tau^n_{l+1}\bigr) -
r^n_l,
\end{eqnarray*}
where in the last line we have used the fact that $\tilde{\mathbf
{X}}^n$ agrees
with $\mathbf{X}^n$ on $\mathcal{P}_n$, as well as the recursive
definition of the
rough path scheme $Y^n$. Hence, we have that
%
\[
\bigl|Z^{(l)}_{l+1} - Z^{(l+1)}_{l+1}\bigr| \lesssim1
\wedge \bigl(\bigl\vert \!\bigl\vert \!\bigl\vert\tilde{\mathbf{X}} ^n \bigr\vert \!\bigr\vert \!\bigr\vert_\gamma^3 + C_n^3 \bigr)
\Delta_n^{3\gamma} \lesssim \bigl(1\wedge C_n^3
\bigr) \Delta _n^{3\gamma}.
\]
It follows from \eqref{eZ1} and \eqref{eZ2} that
%
\[
\bigl|\tilde{Y}^n \bigl(\tau^n_k \bigr) -
Y^n \bigl(\tau^n_k \bigr)\bigr| \lesssim \bigl(1
\wedge C_n^4 \bigr) N_n \Delta_n^{3\gamma}
\lesssim \bigl(1\wedge C_n^4 \bigr)
\Delta_n^{3\gamma-1},
\]
where in the last inequality we have used the assumption $\sup_n N_n
\Delta_n < \infty$. This completes the proof.
\end{pf}

\begin{pf*}{Proof of Theorem \protect\ref{thmmod}}
All that is required is to show that $\tilde{Y}^n = \Phi(\tilde
{\mathbf{X}}^n)$
solves \eqref{eXZ} where $\tilde{\mathbf{X}}^n = (\tilde
{X}^n,\tilde{\mathbb{X}}^n) \in
\mathcal{C}^\gamma([0,T]; \mathbb{R}^d)$ is derived in Lem\-ma~\ref
{thmbxtilde}. By definition and
by construction of $\tilde{\mathbf{X}}^n$ we have that
\[
\tilde{Y}^n(t) = \tilde{Y}^n(s) + V \bigl(
\tilde{Y}^n(s) \bigr)\tilde{X}^n(s,t) + \mathbb{V} \bigl(
\tilde{Y}^n(s) \bigr) \dvtx \tilde{\mathbb{X}}^n(s,t) +
o\bigl(|t-s|\bigr),
\]
where $\tilde{X}^n$ is a piecewise smooth path (obtained from the
signature realizing $g$) and
%
\[
\tilde{\mathbb{X}}^n(s,t) = \int_s^t
\tilde{X}^n(s,r) \otimes d\tilde{X}^n(r) +
\tilde{Z}^n(t) - \tilde{Z}^n(s),
\]
where the integral is of Riemann--Stieltjes type and where $\tilde{Z}^n$
is constructed by concatenating the increments $z(s,t)$, in particular\vspace*{1pt}
$\tilde{Z}^n$ is piecewise Lipschitz. By \cite{friz10}, Theorem~12.14,  it
follows that $\tilde{Y}^n$ satisfies \eqref{eXZ}. Note that
\cite{friz10}, Theorem~12.14,  is basically Lemma~\ref{lemrpsde} but under
the assumption that the driving path is piecewise smooth rather than a
semi-martingale.
\end{pf*}

\subsection{Discrete Kolmogorov criterion}
In Section~\ref{sconvergence}, we will employ the standard method of
lifting weak convergence in the sup-norm topology to weak convergence
in some $\gamma$-H\"older topology, using a tightness condition. In the
\textit{continuous time} setting (which we cannot use), the
Kolmogorov--Lamperti criterion \cite{lamperti62,friz10} is the usual
method for checking this tightness condition. The following is a slight
modification of a version of the criterion found in Corollary A12 \cite{friz10}.

\begin{lemma}\label{thmctskolm}
Let $\mathbf{X}^n = (X^n,\mathbb{X}^n)$ define a sequence of rough
paths. Suppose that
%
\begin{equation}
\label{ectskolm}
\bigl( \mathbf{E}\bigl|X^n(s,t)\bigr|^q
\bigr)^{1/q} \lesssim|t-s|^{\alpha} \quad\mbox{and}\quad \bigl( \mathbf{E}|
\mathbb{X}^n(s,t)|^{q/2} \bigr)^{2/q} \lesssim
|t-s|^{2\alpha}
\end{equation}
for each $s,t \in[0,T]$, uniformly in $n\geq1$. Then
%
\[
\sup_{n\geq1 }\mathbf{E}\bigl\vert \!\bigl\vert \!\bigl\vert\mathbf{X}^n
\bigr\vert \!\bigr\vert \!\bigr\vert_\gamma^q < \infty
\]
for any $\gamma\in(0,\alpha-q^{-1})$. In particular, we have that
%
\begin{equation}
\label{ectstightness}
\sup_{n\geq1} \mathbf{P} \bigl( \bigl\vert \!\bigl\vert \!\bigl\vert
\mathbf{X}^n \bigr\vert \!\bigr\vert \!\bigr\vert_\gamma> M \bigr) \to0\qquad \mbox{as $M \to\infty$.}
\end{equation}
And moreover $\{ \mathbf{X}^n\}_{n\geq1}$ is tight in the $\rho
_\gamma$
topology for every $\gamma\in(0,\alpha-q^{-1})$.
\end{lemma}
\begin{pf}
In the case of geometric rough paths [where $\mathbf{X}^n$ takes
valued in
$G^2(\mathbb{R}^d)$], the result is simply Corollary A12 of \cite{friz10}.
To extend the result to general rough paths, one simply applies the
Garcia--Rodemich--Rumsey interpolation result to the components $X$ and
$\mathbb{X}$ individually. This argument can be found in \cite{gubinelli04},
Corollary~4.
\end{pf}
Obviously, this result cannot be used directly on rough step functions,
since step functions have no hope of satisfying the Kolmogorov
estimates. Fortunately, a discrete version of the above result turns
out to be equally as useful. We define the \textit{discrete tightness
condition} as
%
\begin{equation}
\label{edtight}
\sup_{n\geq1} \mathbf{P} \bigl( \bigl\vert \!\bigl\vert \!\bigl\vert
\mathbf{X}^n \bigr\vert \!\bigr\vert \!\bigr\vert_{\gamma,n} > M \bigr) \to0
\qquad \mbox{as $M \to\infty$.}
\end{equation}
This essentially says that the rough step functions are ``H\"older
continuous,'' provided we do not look at them too closely (i.e.,
near the jumps). We will now show that the discrete tightness criterion
can likewise be checked using a discrete version of the continuous
Kolmogorov criterion. In particular, we need only check the estimate on
the partition $\mathcal{P}_n$.

\begin{lemma}\label{thmkolm}
Suppose that
%
\[
\bigl( \mathbf{E}\bigl|X^n \bigl(\tau^n_j,
\tau^n_k \bigr)\bigr|^q \bigr)^{1/q}
\lesssim \bigl|\tau^n_j - \tau^n_k\bigr|^{\alpha}
\quad \mbox{and}\quad \bigl( \mathbf {E}\bigl|\mathbb{X}^n \bigl(\tau
^n_j,\tau^n_k
\bigr)\bigr|^{q/2} \bigr)^{2/q} \lesssim\bigl|\tau^n_j
- \tau ^n_k\bigr|^{2\alpha}
\]
for each $\tau^n_j,\tau^n_k\in\mathcal{P}_n$ uniformly in $n\geq
1$, for some
$\alpha\in(0,1/2]$. Then the discrete tightness condition \eqref{edtight} holds for any $\gamma\in(0,\alpha- q^{-1})$.
\end{lemma}

\begin{pf}
The idea behind the proof is to replace $\mathbf{X}^n$ with the
genuine rough
path $\tilde{\mathbf{X}}^n$ constructed in Lemma~\ref{thmbxtilde},
which, as you
recall, agrees with $\mathbf{X}^n$ on $\mathcal{P}_n$. Since
%
\[
\bigl\vert \!\bigl\vert \!\bigl\vert\mathbf{X}^n \bigr\vert \!\bigr\vert \!\bigr\vert_{\gamma
,n} =
\bigl\vert \!\bigl\vert \!\bigl\vert\tilde{\mathbf{X}}^n \bigr\vert \!\bigr\vert \!\bigr\vert_{\gamma,n}
\leq\bigl\vert \!\bigl\vert \!\bigl\vert\tilde{\mathbf{X}}^n \bigr\vert \!\bigr\vert \!\bigr\vert_{\gamma},
\]
to prove the discrete tightness condition \eqref{edtight} it is
sufficient to check the H\"older estimate \eqref{ectskolm} for the
process $\tilde{\mathbf{X}}^n$ and apply Lemma~\ref{thmctskolm}.
Hence, we need
only verify that
%
\begin{equation}
\label{etildekolm}
\mathbf{E}\bigl|\tilde{X}^n(s,t)\bigr|^{q}
\lesssim|s-t|^{q\alpha} \quad\mbox{and}\quad \mathbf{E}\bigl|\tilde{\mathbb{X}}^n(s,t)\bigr|^{q/2}
\lesssim |t-s|^{q\alpha},
\end{equation}
holds for each $s,t \in[0,T]$, uniformly in $n\geq1 $.

Assume without loss of generality that $s,t \in[0,T]$ and $\tau
^n_{j-1} < s < \tau^n_{j}$ and $\tau^n_k \leq t < \tau^n_{k+1}$ (note
that the case $s,t \in[\tau_j,\tau_{j+1}]$ is essentially a
sub-argument of the arguments below). From Chen's relations, we know that
%
\begin{equation}
\label{efeb9}
\mathbf{E}\bigl|\tilde{X}^n(s,t)\bigr|^q \lesssim
\mathbf{E}\bigl|\tilde {X}^n \bigl(s,\tau^n_j
\bigr)\bigr|^q + \mathbf{E} \bigl|X^n_{j,k}\bigr|^q
+ \mathbf{E}\bigl|\tilde{X}^n \bigl(\tau^n_k,t
\bigr)\bigr|^q.
\end{equation}
But from \eqref{efeb9a}, we see that
\begin{eqnarray*}
\mathbf{E}\bigl|\tilde{X}\bigl(s,\tau^n_j
\bigr)\bigr|^q &\lesssim &  \bigl(\mathbf {E}\bigl|X^n_{j-1,j}\bigr|^q
+ \mathbf{E} \bigl|\mathbb{X}^n_{j-1,j}\bigr|^{q/2} \bigr)
\biggl(\frac{\tau^n_j-s}{\tau
^n_{j}-\tau
^n_{j-1}} \biggr)^q
\\
&\lesssim & \bigl(\tau^n_j-\tau^n_{j-1}
\bigr)^{q\alpha} \biggl(\frac{\tau^n_j-s}{\tau^n_{j}-\tau^n_{j-1}} \biggr)^q
\\
&=& \bigl(\tau ^n_j - s\bigr)^{q\alpha} \biggl(
\frac{\tau^n_j-s}{\tau^n_j-\tau^n_{j-1}} \biggr)^{q-q\alpha} \leq\bigl(\tau^n_j
- s\bigr)^{q\alpha} \leq(t - s)^{q\alpha}.
\end{eqnarray*}
By assumption, we have that
%
\[
\mathbf{E}\bigl|X^n_{j,k}\bigr|^q \lesssim \bigl(
\tau^n_k-\tau^n_j
\bigr)^{q\alpha} \lesssim (t-s)^{q\alpha}.
\]
The remaining term in \eqref{efeb9} can be bounded similarly. By
Chen's relations (and H\"older's inequality), we also have that
%
\begin{eqnarray}
&&\mathbf{E}\bigl|\tilde{\mathbb{X}}^n(s,t)\bigr|^{q/2} \nonumber\\
&& \label{efeb10}\qquad \lesssim  \mathbf {E}\bigl|\tilde{\mathbb{X}}^n\bigl(s,\tau
^n_j\bigr)\bigr|^{q/2} + \mathbf{E}\bigl|
\mathbb{X}^n_{j,k}\bigr|^{q/2} + \mathbf {E}\bigl|\tilde{
\mathbb{X}}^n\bigl(\tau ^n_k,t
\bigr)\bigr|^{q/2}
\\
\nonumber
&&\qquad\quad{}+ \bigl(\mathbf{E}\bigl|X^n_{j,k}\bigr|^{q}\mathbf{E}\bigr|
\tilde {X}^n\bigl(\tau ^n_k,t
\bigr)\bigr|^q\bigr)^{1/2} + \bigl(\mathbf{E}\bigl|
\tilde{X}^n\bigl(s,\tau^n_j
\bigr)\bigr|^q\mathbf {E}|\tilde{X}^n\bigl(\tau
^n_j,t\bigr)|^q\bigr)^{1/2}.\hspace*{-24pt}
\end{eqnarray}
But from \eqref{efeb10a} we have that
\begin{eqnarray*}
\mathbf{E}\bigl|\tilde{\mathbb{X}}^n\bigl(s,\tau^n_j
\bigr)\bigr|^{q/2} &\lesssim & \biggl(\frac{t-s}{\tau
^n_{j+1}-\tau^n_j} \biggr)^q
\bigl( \mathbf{E}\bigl|X^n_{j,j+1}\bigr|^q + \mathbf{E}\bigl|
\mathbb{X} ^n_{j,j+1}\bigr|^{q/2} \bigr)
\\
&&{}+ \biggl(\frac{t-s}{\tau^n_{j+1}-\tau
^n_j} \biggr)^{q/2} \mathbf{E}\bigl|
\mathbb{X}^n_{j,j+1}\bigr|^{q/2}.
\end{eqnarray*}
As above, it follows that
%
\[
\mathbf{E}\bigl|\tilde{\mathbb{X}}^n \bigl(s,\tau^n_j
\bigr)\bigr|^{q/2} \lesssim \bigl(\tau ^n_j-s
\bigr)^{q\alpha} \leq (t-s)^{q\alpha}.
\]
The other terms in \eqref{efeb10} can be bounded similarly. This
completes the proof.
\end{pf}

%

\begin{rmk}
The discrete criterion differs from the continuous case in the
assumption $\alpha\leq1/2$, which was not required in the continuous
case. However, this assumption only becomes a restriction when the
diffusion approximation is driven by a path with H\"older exponent
$\gamma> 1/2$. Of course, one can always resolve the problem by
treating the path as having the weaker H\"older exponent. On the other
hand, in these higher regularity situations the iterated integrals
become unnecessary and a much simpler theory of Young integration (with
much weaker assumptions) would suffice.
\end{rmk}

\section{Convergence of rough path schemes}\label{sconvergence}

We can now prove the main result of the article.
%
\begin{thm}\label{thmweak}
Suppose that $\mathbf{X}^n \stackrel{\mathit{f.d.d.}}{\to}\mathbf{X}$
and that $\mathbf{X}^n$ satisfies the discrete
tightness condition \eqref{edtight} for some $\gamma\in(1/3,1/2]$.
Then $Y^n = \Phi^n(\mathbf{X}^n) \Rightarrow\Phi(\mathbf{X})$ in
the sup-norm topology.
\end{thm}

\begin{pf}
First,\vspace*{1pt} let $\tilde{\mathbf{X}}^n$ be the $\gamma$-H\"older rough path
constructed in Lem\-ma~\ref{thmbxtilde} and let $\tilde{Y}^n = \Phi
(\tilde{\mathbf{X}}^n)$. To prove the theorem, it is sufficient to
first show
that $\llVert  Y^n-\tilde{Y}^n \rrVert _\infty\to0$ in probability
and second
show that $\tilde{Y}^n \Rightarrow Y$ in the sup-norm topology, hence we
will proceed as such.

As usual, we use the shorthand $C_n = \vert \!\vert \!\vert\mathbf
{X}^n \vert \! \vert \! \vert_{\gamma,n}$. From
Lemma~\ref{thmytilde}, it follows that
%
\[
\bigl\llVert Y^n - \tilde{Y}^n \bigr\rrVert
_\infty\lesssim \bigl(1\wedge C_n^4 \bigr)
\Delta _n^{3\gamma-1}.
\]
Hence,
\begin{eqnarray*}
\mathbf{P} \bigl( \bigl\llVert Y^n - \tilde{Y}^n \bigr
\rrVert _\infty> \delta \bigr) &\leq &\mathbf{P} \bigl( C \bigl(1\wedge
C_n^4\bigr)\Delta_n^{3\gamma-1} > \delta
\bigr)
\\
&=& \mathbf{P} \bigl( 1\wedge\bigl\vert \!\bigl\vert \!\bigl\vert\mathbf{X}^n \bigr\vert \!\bigr\vert \!\bigr\vert_{\gamma,n} > \bigl(C^{-1}\delta\Delta
_n^{1-3\gamma} \bigr)^{1/4} \bigr).
\end{eqnarray*}
But since $\Delta_n^{1-3\gamma} \to\infty$ as $n\to\infty$, we see
that for any arbitrarily large $M>0$
\begin{eqnarray*}
\lim_{n\to\infty} \mathbf{P} \bigl( \bigl\llVert Y^n -
\tilde{Y}^n \bigr\rrVert _\infty> \delta \bigr) &\leq &\limsup
_{n\to\infty} \mathbf{P} \bigl( \bigl\vert \!\bigl\vert \!\bigl\vert\mathbf{X}^n
\bigr\vert \!\bigr\vert \!\bigr\vert_{\gamma,n} > C^{-1}\delta^{1/4}
\Delta_n^{(1-3\gamma)/4} \bigr)
\\
&\leq & \limsup_{n\to\infty}\mathbf{P} \bigl( \bigl\vert \!\bigl\vert \!\bigl\vert
\mathbf{X}^n \bigr\vert \!\bigr\vert \!\bigr\vert_{\gamma,n} > M \bigr).
\end{eqnarray*}
Finally, by taking $M\to\infty$, it follows from the discrete tightness
condition that
%
\[
\lim_{n\to\infty} \mathbf{P} \bigl( \bigl\llVert Y^n -
\tilde{Y}^n \bigr\rrVert _\infty> \delta \bigr) = 0.
\]
Now we prove that $\tilde{Y}^n \Rightarrow Y$ in the sup-norm topology.
Due to the continuity of the map $\Phi$, as stated in Lemma~\ref
{thmitomap}, it is sufficient to prove that $\tilde{\mathbf{X}}^n
\Rightarrow\mathbf{X}
$ in the $\rho_\gamma$ topology. It is therefore sufficient to first
show that $\tilde{\mathbf{X}}^n\stackrel{\mathrm{f.d.d.}}{\to}\mathbf
{X}$ and second that $\{\tilde{\mathbf{X}}^n\}_{n\geq
1}$ is tight in the $\rho_\gamma$ topology.

First, due to the regularity of $\tilde{\mathbf{X}}^n$ between mesh
points, it is
easy to see that $\llVert  \mathbf{X}^n(t) - \tilde{\mathbf{X}}^n(t)
\rrVert _\infty\lesssim
(1\vee C_n^2) \Delta_n^{\kappa}$ for some $\kappa>0$. Hence, by an
argument similar to that found at the start of the proof, it follows
from the discrete tightness condition that $\llVert  \mathbf{X}^n -
\tilde{\mathbf{X}} ^n \rrVert _\infty\to0$ in probability. And
since by assumption $\mathbf{X}^n \stackrel{\mathrm{f.d.d.}}{\to}
\mathbf{X}$ it follows that $\tilde{\mathbf{X}}^n \stackrel{\mathrm
{f.d.d.}}{\to}\mathbf{X}$. We will now move onto the
tightness argument. From Lemma~\ref{thmbxtilde}, we have the estimate
$\vert \!\vert \!\vert\tilde{\mathbf{X}}^n \vert  \!\vert \! \vert
_\gamma\lesssim\vert \!\vert \!\vert\mathbf{X}^n \vert  \!\vert
\!\vert_{\gamma,n}$. It
follows that
%
\[
\mathbf{P} \bigl( \bigl\vert \!\bigl\vert \!\bigl\vert\tilde{\mathbf{X}}^n \bigr\vert \!\bigr\vert \!\bigr\vert_{\gamma} > M \bigr)
\leq\mathbf{P} \bigl( C \bigl\vert \!\bigl\vert \!\bigl\vert
\mathbf{X}^n \bigr\vert \!\bigr\vert \!\bigr\vert_{\gamma
,n} > M \bigr),
\]
and the tightness of $\tilde{\mathbf{X}}^n$ in the $\rho_\gamma$
topology follows
from the discrete tightness condition.
This completes the proof.
\end{pf}

We can now prove the theorems introduced in Section~\ref{sresults}.
They are both immediate corollaries.

\begin{pf*}{Proof of Theorem \protect\ref{thmmain}}
$\!$By the discrete Kolmogorov criterion (Lem\-ma~\ref{thmkolm}), we obtain
the discrete tightness criterion, and hence can apply Theorem~\ref
{thmweak}. To identify the limit $Y = \Phi(\mathbf{X})$, we simply apply
Lemma~\ref{lemrpsde}.
\end{pf*}

\begin{rmk}\label{rmkqlarge}\label{rmkdrift2} To prove the result
with the relaxed assumption described in Remark~\ref{rmkqsmall}, one
simply replaces Lemma~\ref{thmitomap} with the sharper version \cite{friz10}, Theorem~12.10,  and the remaining argument is identical. To
prove the result with additional drift, as described in Remark~\ref
{rmkdrift1}, we again replace Lemma~\ref{thmitomap} with \cite{friz10}, Theorem~12.10, but now we must use the $(p,q)$-rough path
$(\tilde{\mathbf{X}}^n,t)$ with $(p,q) = (2-\kappa,1)$ and
arbitrarily small
$\kappa> 0$. The remaining argument is identical.
\end{rmk}

\begin{pf*}{Proof of Meta Theorem \protect\ref{thmmeta}}
The proof is completely identical to the proof above, but we still need
to ``interpret'' the limit $Y = \Phi(\mathbf{X})$. This is a fairly
nonrigorous statement and, therefore, has a fairly nonrigorous proof.
We are merely sketching an idea that would apply more rigorously in
concrete situations.

By definition, the limit $Y$ solves the RDE
%
\begin{equation}
\label{eRDEfinal}
Y(t) = \int_0^t V \bigl(Y(s)
\bigr) \,d\mathbf{X}(s).
\end{equation}
It is a general \textit{heuristic} that if $\mathbb{X}$ is constructed using
some known construction then the integral in \eqref{eRDEfinal} is
constructed similarly. For instance, suppose there is some ``method of
integration,'' which is a bilinear operator
%
\[
\int_0^t A \star dB \stackrel{{\mathrm{def}}}
{=}\mathcal{I}(A,B) (t)
\]
for two continuous paths $A,B \in C([0,T];\mathbb{R})$ satisfying the
obvious condition
%
\[
\mathcal{I}(1,A) (t) = A(t).
\]
Now suppose that $\mathbb{X}$ is defined by
%
\[
\mathbb{X}^{\alpha\beta}(t) = \mathcal{I} \bigl(X^\alpha,X^\beta
\bigr) (t),
\]
for each $\alpha,\beta= 1,\ldots,d$. Then, using the theory of
controlled rough paths \cite{gubinelli04,gubinelli10a}, it can be shown
that $Y$ solves \eqref{eRDEfinal} if and only if $Y$ is a fixed point
of the equation
%
\[
Y(t) = Y(0) + \int_0^t V \bigl(Y(s) \bigr)
\star dX(s).
\]
The assumptions on $\mathcal{I}$ are generic enough to include
virtually any
reasonable construction of an integration map (for integrators with H\"
older exponent $\gamma> 1/3$).
\end{pf*}

\section*{Acknowledgements}
David Kelly would like to thank I.~Melbourne for introducing the motivating
problem and P.~Friz for several constructive comments.

\printaddresses
\end{document}